\documentclass{amsart}
\usepackage{graphicx}
\usepackage[all]{xy}
\usepackage{latexsym}
\usepackage{amssymb}
\usepackage{amsmath}
\usepackage{color}

\xyoption{arc}

 \newfam\cyrfam

  \font\tencyr=wncyr10

  \font\sevencyr=wncyr7

  \font\fivecyr=wncyr5

  \textfont\cyrfam=\tencyr \scriptfont\cyrfam=\sevencyr

    \scriptscriptfont\cyrfam=\fivecyr

  \newfam\cyifam

  \font\tencyi=wncyi10

  \font\sevencyi=wncyi7

  \font\fivecyi=wncyi5

  \textfont\cyifam=\tencyi \scriptfont\cyifam=\sevencyi

    \scriptscriptfont\cyifam=\fivecyi

\def\id{{\mbox{1 \hskip -7pt 1}}}
\newcommand{\sgn}{{\mathit s  \mathit g\mathit  n}}
 \newcommand{\lon}{\longrightarrow}
 \newcommand{\bu}{\bullet}
 
 \newcommand{\rar}{\rightarrow}

\newcommand{\p}{{\partial}}

\newcommand{\tw}{\mathsf{tw}}
\newcommand{\RGra}{\mathcal{R} \mathcal{G} ra}

\newcommand{\Grav}{\mathcal{G} rav}

\newcommand{\ST}{\mathcal{ST}}

 \newcommand{\Z}{{\mathbb Z}}
 \newcommand{\bS}{{\mathbb S}}

 \newcommand{\R}{{\mathbb R}}
 
 \newcommand{\K}{{\mathbb K}}

 \newcommand{\ot}{\otimes}

\newcommand{\sE}{{\mathsf E}}

\newcommand{\sR}{{\mathsf R}}

\newcommand{\Def}{\mathsf{Def}}

\newcommand{\GC}{\mathsf{GC}}

 \newcommand{\Beq}{\begin{equation}}
 \newcommand{\Eeq}{\end{equation}}
 \newcommand{\Beqr}{\begin{eqnarray}}
 \newcommand{\Eeqr}{\end{eqnarray}}
 \newcommand{\Beqrn}{\begin{eqnarray*}}
 \newcommand{\Eeqrn}{\end{eqnarray*}}
 \newcommand{\Ba}{\begin{array}}
 \newcommand{\Ea}{\end{array}}
 \newcommand{\Bi}{\begin{itemize}}
 \newcommand{\Ei}{\end{itemize}}
 \newcommand{\Bc}{\begin{center}}
 \newcommand{\Ec}{\end{center}}
 \newcommand{\fo}{{\mathfrak o}}

\newcommand{\fs}{{\mathfrak s}}

 \newcommand{\cA}{{\mathcal A}}
 
 \newcommand{\cC}{{\mathcal C}}
 
 \newcommand{\cE}{{\mathcal E}}
 \newcommand{\cF}{{\mathcal F}}
 \newcommand{\cG}{{\mathcal G}}

 \newcommand{\cM}{{\mathcal M}}
 \newcommand{\cN}{{\mathcal N}}
 \newcommand{\cP}{{\mathcal P}}
 
 \newcommand{\cR}{{\mathcal R}}
 
 \newcommand{\cT}{{\mathcal T}}

 \newcommand{\al}{\alpha}
 \newcommand{\be}{\beta}
 \newcommand{\ga}{\gamma}
 
 \newcommand{\Ga}{\Gamma}

 \newcommand{\om}{\omega}

 \newcommand{\Ker}{{\mathsf K \mathsf e \mathsf r}\, }
 
 \newcommand{\Hom}{{\mathrm H\mathrm o\mathrm m}}
 
 \newcommand{\sip}{\smallskip}
 \newcommand{\bip}{\bigskip}
 \newcommand{\mip}{\vspace{2.5mm}}
\newcommand{\LB}{\mathcal{L}\mathit{ieb}}
\newcommand{\LoB}{\mathcal{L}\mathit{ieb}^\diamond}

\newcommand{\LBd}{\mathcal{L}\mathit{ieb}_{d}}

\newcommand{\LoBd}{\mathcal{L}\mathit{ieb}_{d}^\diamond}

\newcommand{\HoLoBd}{\mathcal{H}\mathit{olieb}_{d}^\diamond}

\newcommand{\HoLoB}{\mathcal{H}\mathit{olieb}^\diamond}

\newcommand{\Lie}{\mathcal{L} \mathit{ie}}

\theoremstyle{plain}
\swapnumbers

\newtheorem{prop-def}[theorem]{Proposition-definition}

\newtheorem{f-theorem}{Formality Theorem}[section]
\newtheorem{main-theorem}{Main~Theorem}[section]
\newtheorem{section-theorem}{Theorem}[section]

\theoremstyle{definition}

\begin{document}

 \sloppy

 \newenvironment{proo}{\begin{trivlist} \item{\sc {Proof.}}}
  {\hfill $\square$ \end{trivlist}}

\long\def\symbolfootnote[#1]#2{\begingroup%
\def\thefootnote{\fnsymbol{footnote}}\footnote[#1]{#2}\endgroup}

\title{From gravity to string topology}

\author{Sergei~A. Merkulov}
\address{Sergei~Merkulov:  Department of Mathematics, Luxembourg University,  Grand Duchy of Luxembourg }
\email{sergei.merkulov@uni.lu}

\begin{abstract}
The chain gravity properad %$\cC h\Grav_{3-d}$
introduced in \cite{Me2} acts on the cyclic Hochschild of any cyclic $A_\infty$ algebra equipped with a scalar product of degree $-d$. In particular, it acts on the cyclic Hochschild complex of any  Poincare duality algebra of degree $d$, and that action factors through a quotient dg properad $\ST_{3-d}$ of ribbon graphs which is in focus of this paper. We show that its cohomology properad $H^\bu(\ST_{3-d})$ is highly non-trivial and that it acts canonically on the reduced equivariant homology $\bar{H}_\bu^{S^1}(LM)$ of the loop space of any  simply connected $d$-dimensional closed manifold $M$. By its very construction, the string topology properad $H^\bu(\ST_{3-d})$ comes equipped with a morphism from the
  gravity properad
 $
 \Grav_{3-d}%:=H^\bu(\cC h\Grav_{3-d})=
 %\left\{\prod_{g,m,n}H_c^{\bu+(3-d)(2g-2+m+n)}(\cM_{g,m+n}\times \R^m)\right\},
  $
  which is fully determined by the compactly supported cohomology of the moduli spaces $\cM_{g,n}$  of stable algebraic curves of genus $g$ with marked points. This result gives rise to new universal operations in string topology as well as reproduces in a unified way several known constructions:  we show that (i) $H^\bu(\ST_{3-d})$ is also a properad under the properad of involutive Lie bialgebras $\LoB_{3-d}$ whose induced (via  $H^\bu(\ST_{3-d})$) action on $\bar{H}_\bu^{S^1}(LM)$ agrees precisely with the famous purely geometric construction of M.\ Chas and D.\ Sullivan \cite{CS1,CS2}, (ii)  $H^\bu(\ST_{3-d})$ is a properad under the properad of {\em homotopy}\, involutive Lie bialgebras $\HoLoB_{2-d}$ which controls (via $H^\bu(\ST_{3-d})$) four universal  string topology operations introduced in \cite{Me1},
  (iii) E.\ Getzler's gravity operad
  %$$
  %grav_d:=\{H_c^{\bu+(3-d)(n-1)-1}(\cM_{0,1+n})\}_{n\geq 2}
  %$$
   {\em injects}\, into  $H^\bu(\ST_{3-d})$
    implying a purely algebraic counterpart of the geometric construction of C.\ Westerland \cite{We} establishing an action of the gravity operad on $\bar{H}_\bu^{S^1}(LM)$.
    %DEGREES
    %$H_0$ corresponds to $H_c^{top=\dim}$, we have \dim=2(1+n) - 3; so the degree of the $n$th %operation is $\bu=2+2n - 3 -3n+3+dn-d +1=-n+dn-d=d(n-1)-n$ --- in full agreement with [We]

\sip

%\sip
%\noindent {\sc Mathematics Subject Classifications} (2000). 17B37, 16W30, 51M20.

%\noindent {\sc Key words}. Strongly homotopy bialgebras, Hopf algebras,  Lie bialgebras,
% deformation quantization, configuration spaces.
\end{abstract}

 \maketitle
 \markboth{}{}

{\small
{\small
\tableofcontents
}
}

{\Large
\section{\bf Introduction}
}

This paper aims to show new results in a part of string topology which deals with universal operations on the
 (reduced) equivariant equivariant homology ${H}_\bu^{S^1}(LM)$ of the free loop space $LM$ of a connected and simply-connected closed $d$-dimensional manifold $M$. A large family of operations  on ${H}_\bu^{S^1}(LM)$ has been found already in the fundamental work \cite{CS1} by M.\ Chas and D.\ Sullivan; later it was shown by C.\ Westerland in \cite{We} that these operations provide us with a representation of the gravity operad $grav_{3-d}$ introduced by E.\ Getzler in \cite{Ge}.

\sip

M.\ Chas and D.\ Sullivan found a geometric construction \cite{CS2} of the involutive Lie bialgebra structure  on the reduced equivariant equivariant homology $\bar{H}_\bu^{S^1}(LM)$ of the free loop space a closed $d$-dimensional manifold; in the properadic language adopted in this paper, their construction gives us a representation of the properad $\LoB_{3-d}$ on $\bar{H}_\bu^{S^1}(LM)$ (see \S2 below for details). In the case $d=2$ their construction reproduces the famous Goldman-Turaev Lie bialgebra structure \cite{Go,Tu}  on the reduced equivariant homology of free loop spaces of Riemann surfaces; the latter Lie bialgebra structure can be lifted to the full homology (with constant loops including) and leads us to the beautiful formality theory of the geometric and purely algebraic counterparts of this structure which has been developed by A.\ Alekseev, N.\ Kawazumi, Y.\ Kuno,
G.\ Massuyeau and F.\ Naef in \cite{AKKN1, AKKN2,AN, Mas}.

\sip

A nice purely algebraic construction of the $\LoB_{3-d}$-algebra structure on $\bar{H}_\bu^{S^1}(LM)$ for simply connected closed manifolds $M$ has been found by X.\ Chen, F.\ Eshmatov and W.\ L. Gan in \cite{CEG}. It was further studied by F.\ Naef and T.\ Willwacher in \cite{NW} who showed, in particular, that this purely algebraic construction gives us the same operations as the geometric one; moreover, they extended it to not necessarily simply connected manifolds $M$ using a partition function $Z_M$ introduced earlier in \cite{CW}.
Props of directed ribbon graphs of various types have been introduced and studied in \cite{TZ} with the purpose to get new string topology operations from higher homotopy products on  Hochschild (co)homologies of cyclic $\cA ss_\infty$ algebras.

\sip

In this paper we study string topology applications of the gravity chain properad $\cC h\Grav_d$ introduced by the author in \cite{Me2}; it is generated by ribbon graphs with vertices of two types, white vertices which are labelled and black vertices which are unlabelled and are assigned the cohomological degree $d$; the differential in $\cC h\Grav_d$ creates a new black vertex (see \S 2 below for full details). As explained in \S 3, this dg properad acts (almost by its very construction) on the Hochschild complex $Cyc(A^*[-1])$ of any cyclic $\cA ss_\infty$ algebra $A$ equipped with a degree $-d$ scalar product. The properad $\cC h\Grav_d$  has interesting cohomology,
\Beq\label{1: def of Grav_d}
\Grav_d:= H^\bu(\cC h\Grav_d)\simeq \left\{\prod_{g\geq 0\atop 2g+m+n\geq 3} H^{\bu-m +d(2g-2+m+n)}_c(\cM_{g,m+n})\ot \sgn_m\right\}_{m\geq 1, n\geq 0}
\Eeq
which is called the {\em gravity properad}\, and which is fully determined by the compactly supported cohomology of moduli spaces $\cM_{g,m+n}$ of algebraic curves with $m+n$ marked points ($m$ of them are called {\em out}-points and $n$ of them are called {\em in}\, points). Moreover, the gravity properad $\Grav_d$ contains E.\ Getzler's gravity operad
  \Beq\label{1: grav_d operad}
  grav_d:=\{H_c^{\bu+(3-d)(n-1)-1}(\cM_{0,1+n})\}_{n\geq 2}
  \Eeq
  as a sub-properad, and the latter is known to play a role in string topology \cite{We}. Our main purpose in this paper is to study string topology applications  of its properadic extension $\cG rav_d$, more precisely, of its chain version  $\cC h\Grav_d$ which is more informative and useful in this respect.

\sip

 Given a Poincar\'e duality algebra $A$ in degree $d$, there is a canonical representation of the chain gravity properad (see \S 3 for full details),
$$
\rho_A: \cC h\Grav_{3-d} \lon \cE nd_{ \overline{Cyc}(\bar{A}^*[-1])}.
$$
in the reduced cyclic Hochschild of $A$,
$$
\overline{Cyc}(\bar{A}^*[-1])\simeq \bigoplus_{k\geq 1} \left(\ot^k (\bar{A}^*[-1])\right)_{\Z^*_k}.
$$
Our interest in this class of representations stems from the fact that if $A$ is a Poincar\'e duality model of a connected  $d$-dimensional closed manifold $M$, then there is a linear map
$$
\bar{H}_\bu^{S^1}(LM) \lon H^{\bu}( \overline{Cyc}(\bar{A}^*[-1])
$$
from the reduced equivariant equivariant homology $\bar{H}_\bu^{S^1}(LM)$ of the free loop space $LM$ of $M$ into the cyclic Hochschild cohomology. Moreover, if $M$ is simply connected, this map is an isomorphism (moreover, in this case a Poincar\'e duality model $A$ for $M$ always exists \cite{LS}) so that the  gravity properad $\Grav_d$ acts
canonically on $\bar{H}_\bu^{S^1}(LM)$.

 \sip

 We study in \S 4 the joint kernel of the above class of representation maps $\rho_A$ for all possible $A$
 which leads us to an observation that any such representation $\rho_A$ factors through the canonical epimorphism
 $$
 \cC h\Grav_d \lon \ST_d
 $$
to a {\em chain string topology properad}\, $\ST_d$ which is defined as the quotient of $Ch\Grav_d$ by the differential closure of the ideal generated by ribbon graphs with at least one black vertex of valency $\geq 4$ or with at least one boundary made of black vertices only.
The cohomology properad $H^\bu(\ST_{3-d})$ acts therefore on $\bar{H}_\bu^{S^1}(LM)$
for any simply-connected closed $d$-dimensional manifold $M$. Our next purpose in this paper is to show that the string topology properad $H^\bu(\ST_{d})$ is very non-trivial, and to establish its connection to known results in string topology.

\sip

In \S 4 we prove that there is
\Bi
\item[(i)] a non-trivial morphism of properads
$$
\LoB_d \lon H^\bu(\ST_d)
$$
which induces, via the above action of $H^\bu(\ST_{3-d})$ on $\bar{H}_\bu^{S^1}(LM)$,
the Chas-Sullivan involutive Lie bialgebra structure;

\item[(ii)] an injection of the gravity operad
$$
grav_d \lon H^\bu(\ST_d)
$$
which gives a purely algebraic counterpart of the geometric construction in \cite{
We} establishing an action of $grav_{3-d}$ of $\bar{H}_\bu^{S^1}(LM)$;

\item[(iii)] a  morphism of dg properads
$$
\HoLoB_{d-1} \lon H^\bu(\ST_d)
$$
which is non-trivial on the quartette of {\em homotopy}\, involutive Lie bialgebra bialgebra operations (\ref{4: four trinary cycles in ST}) found earlier in \cite{Me2} using a very different technique of ribbon {\em hyper}graphs.
\Ei

Thus the string topology properad $H^\bu(\ST_d)$  gives  us both well- and less known universal string topology operations on the equivariant cohomology $\bar{H}_\bu^{S^1}(LM)$ of free loop spaces of simply connected closed manifolds in a unified way. It is an open problem to compute all the cohomology
$H^\bu(\ST_d)$ explicitly, as well as  the cohomology of the deformation complex of the morphism (i) which is a Lie algebra under the cohomology Lie algebra of the mysterious even M.\ Kontsevich graph complex \cite{Ko,W}.

\bip

\subsection{Notation}
 The set $\{1,2, \ldots, n\}$ is abbreviated to $[n]$;  its group of automorphisms is
denoted by $\bS_n$. The cyclic subgroup of $\bS_n$ generated by the permutation $(12...n)$ is denoted by $\Z_n^*$.
%for a non-negative integer $d$, the symbol  $\sgn^d$ stands for the
%one-dimensional representation of $\bS_n$ on which $\bS_n$ acts trivially if $d$ is even
%and by sign if $d$ is odd.
The trivial (resp., sign) one-dimensional representation of
 $\bS_n$ is denoted by $\id_n$ (resp. $\sgn_n$).
 %; we often abbreviate $\sgn_n^d:= \sgn_n^{\ot |d|}$.
 %; we also set, for $d\in \Z$,
 %$$
 %\sgn^d_n:=\left\{ \Ba{ll} \id_n & \mathrm{if}\ d\ \mathrm{is\ even}\\
 %\sgn_n & \mathrm{if}\ d\ \mathrm{is\ odd}.
 %\Ea\right.
 %$$
The cardinality of a finite set $I$ is denoted by $\# I$. A  linear span of a set $I$ over a field $\K$ is denoted by $\K\langle I\rangle$.

\sip

We work throughout in the category of $\Z$-graded vector spaces over a field $\K$
of characteristic zero.
If $V=\oplus_{i\in \Z} V^i$ is a graded vector space, then
$V[k]$ stands for the graded vector space with $V[k]^i:=V^{i+k}$ and
and $s^k$ for the associated isomorphism $V\rar V[k]$; for $v\in V^i$ we set $|v|:=i$. The endomorphism properad\footnote{For a nice introduction into the theory of props and properads we refer to the paper \cite{V} by B.\ Vallette.} of $V$ is denoted by $\cE nd_V$.
%For a pair of graded vector spaces $V_1$ and $V_2$, the symbol $\Hom_i(V_1,V_2)$ stands
%for the space of homogeneous linear maps of degree $i$, and
%$\Hom(V_1,V_2):=\bigoplus_{i\in \Z}\Hom_i(V_1,V_2)$; for example, $s^k\in \Hom_{-k}(V,V[k])$.

\sip

 All our complexes have differential of degree $+1$.
In particular, the cohomology group, $H^\bu(M)$, of a topological space $M$ is non-negatively graded as usual, while its homology group $H_\bu(M)$ is {\em non-positively}\, graded (so that both spaces are dual to each other as $\Z$-graded spaces).
\mip

%{\bf Acknowledgement}.
%It is a pleasure to thank Alexey Kalugin,  Anton Khoroshkin, Sergei Shadrin and especially Thomas %Willwacher  for valuable communications.

\bip

\bip

{\Large
\section{\bf A brief introduction into the (chain) gravity properad}
}

\mip

\subsection{T.\ Willwacher's twisting endofunctor}
For any $d\in \Z$ the operad of degree $d$ shifted Lie algebras is, by definition, the quotient
$$
\Lie_{d}:=\cF ree\langle E\rangle/\langle\cR\rangle,
$$
of the free operad generated by an  $\bS$-module $E=\{E(n)\}_{n\geq 2}$
$$
E(n):=\left\{ \Ba{ll} \sgn_2^{d}\ot \id_1[d-1]=\mbox{span}\left\langle
\Ba{c}\begin{xy}
 <0mm,0.66mm>*{};<0mm,3mm>*{}**@{-},
 <0.39mm,-0.39mm>*{};<2.2mm,-2.2mm>*{}**@{-},
 <-0.35mm,-0.35mm>*{};<-2.2mm,-2.2mm>*{}**@{-},
 <0mm,0mm>*{\bu};<0mm,0mm>*{}**@{},
   %<0mm,0.66mm>*{};<0mm,3.4mm>*{^1}**@{},
   <0.39mm,-0.39mm>*{};<2.9mm,-4mm>*{^{_2}}**@{},
   <-0.35mm,-0.35mm>*{};<-2.8mm,-4mm>*{^{_1}}**@{},
\end{xy}\Ea
=(-1)^{d}
\Ba{c}\begin{xy}
 <0mm,0.66mm>*{};<0mm,3mm>*{}**@{-},
 <0.39mm,-0.39mm>*{};<2.2mm,-2.2mm>*{}**@{-},
 <-0.35mm,-0.35mm>*{};<-2.2mm,-2.2mm>*{}**@{-},
 <0mm,0mm>*{\bu};<0mm,0mm>*{}**@{},
   %<0mm,0.66mm>*{};<0mm,3.4mm>*{^1}**@{},
   <0.39mm,-0.39mm>*{};<2.9mm,-4mm>*{^{_1}}**@{},
   <-0.35mm,-0.35mm>*{};<-2.8mm,-4mm>*{^{_2}}**@{},
\end{xy}\Ea
\right\rangle  & \text{if}\ n=2,\\
0 & \text{otherwise}.
\Ea\right.
$$
 modulo the ideal generated by the following relation\footnote{When representing elements of operads and props as decorated graphs we tacitly assume that all edges and legs are {\em directed}\, along the flow going from the bottom of the graph to the top.

 The action of the element $\sum_{k=1}^{3} (123)^k\in \K[\bS_3]$ on any element $a$ of an $\bS_3$-module is denoted by $\oint_{123} a$.}
$$
%%%%%%%%%%%%%% Lie %%%%%%%%%%%%%%%%%%%%%%%%
\cR: \ \ \ \
%%%%%%%%%%%%%% Lie %%%%%%%%%%%%%%%%%%%%%%%%
\oint_{123} \Ba{c}\resizebox{10mm}{!}{ \begin{xy}
 <0mm,0mm>*{\bu};<0mm,0mm>*{}**@{},
 <0mm,0.69mm>*{};<0mm,3.0mm>*{}**@{-},
 <0.39mm,-0.39mm>*{};<2.4mm,-2.4mm>*{}**@{-},
 <-0.35mm,-0.35mm>*{};<-1.9mm,-1.9mm>*{}**@{-},
 <-2.4mm,-2.4mm>*{\bu};<-2.4mm,-2.4mm>*{}**@{},
 <-2.0mm,-2.8mm>*{};<0mm,-4.9mm>*{}**@{-},
 <-2.8mm,-2.9mm>*{};<-4.7mm,-4.9mm>*{}**@{-},
    <0.39mm,-0.39mm>*{};<3.3mm,-4.0mm>*{^3}**@{},
    <-2.0mm,-2.8mm>*{};<0.5mm,-6.7mm>*{^2}**@{},
    <-2.8mm,-2.9mm>*{};<-5.2mm,-6.7mm>*{^1}**@{},
 \end{xy}}\Ea=0
$$
This operad controls dg Lie algebras $(V,d)$ with the Lie bracket $[\ ,\ ]$
of degree $1-d$. Hence the case $d=1$ corresponds to the ordinary operad of Lie algebras. As usual,  Maurer-Cartan elements of  a $\Lie_d$-algebra $V$ are defined as degree $d$ elements $\ga\in V$ satisfying the equation
\Beq\label{2: MC eqn}
d\ga + \frac{1}{2} [\ga,\ga]=0.
\Eeq

\sip

Given a dg properad $(\cP=\{\cP(m,n)\}_{m,n\geq 0}, \p)$ under $\Lie_d$, that is, a dg properad equipped with a morphism
\Beq\label{2: i from Lie to P}
\imath: \Lie_d \lon \cP,
\Eeq
there is \cite{W}  an associated dg properad $(\tw \cP, \delta)$ which is generated by $\cP$ and one extra   generator with no inputs and one output; we represent this special element pictorially as a corolla
 $
  \Ba{c}\resizebox{1.7mm}{!}{\begin{xy}
 <0mm,0.5mm>*{};<0mm,3mm>*{}**@{-},
 <0mm,0mm>*{\bullet};<0mm,0mm>*{}**@{},
 \end{xy}}\Ea
 $ with black vertex and assign to it the cohomological degree $d$. Generic elements $a$ of $\cP(m,n)$ can also identified with $(m,n)$-corollas
\Beq\label{2: generic elements of cP as (m,n)-corollas}
\Ba{c}\resizebox{15mm}{!}{
 \begin{xy}
 <0mm,0mm>*{\circ};<-8mm,6mm>*{^1}**@{-},
 <0mm,0mm>*{\circ};<-4.5mm,6mm>*{^2}**@{-},
 <0mm,0mm>*{\circ};<0mm,5.5mm>*{\ldots}**@{},
 <0mm,0mm>*{\circ};<3.5mm,5mm>*{}**@{-},
 <0mm,0mm>*{\circ};<8mm,6mm>*{^m}**@{-},
 <0mm,0mm>*{\circ};<-8mm,-6mm>*{_1}**@{-},
 <0mm,0mm>*{\circ};<-4.5mm,-6mm>*{_2}**@{-},
 <0mm,0mm>*{\circ};<0mm,-5.5mm>*{\ldots}**@{},
 <0mm,0mm>*{\circ};<4.5mm,-6mm>*+{}**@{-},
 <0mm,0mm>*{\circ};<8mm,-6mm>*{_n}**@{-},
   \end{xy}}\Ea
\Eeq
whose (say, white) vertex is decorated by $a$; then properadic compositions in $\cP$ can be represented pictorially by gluing out-legs of such decorated corollas to in-legs  of another decorated corollas. Thus the properad $\tw \cP$ is freely generated by white corollas (\ref{2: generic elements of cP as (m,n)-corollas}) and one distinguished $(1,0)$ corolla   \hspace{-0.7mm}$\Ba{c}\resizebox{1.7mm}{!}{\begin{xy}
 <0mm,0.5mm>*{};<0mm,3mm>*{}**@{-},
 <0mm,0mm>*{\bullet};<0mm,0mm>*{}**@{},
 \end{xy}}\Ea$\hspace{-1mm}. The map (\ref{2: i from Lie to P}) sends the generator of $\Lie_d$  into a distinguished degree $1-d$ element of $\cP(1,2)$,
$$
\imath\left(\hspace{-1mm}\Ba{c}\begin{xy}
 <0mm,0.66mm>*{};<0mm,4mm>*{}**@{-},
 <0.39mm,-0.39mm>*{};<2.2mm,-2.2mm>*{}**@{-},
 <-0.35mm,-0.35mm>*{};<-2.2mm,-2.2mm>*{}**@{-},
 <0mm,0mm>*{\bu};<0mm,0mm>*{}**@{},
   %<0mm,0.66mm>*{};<0mm,3.4mm>*{^1}**@{},
   <0.39mm,-0.39mm>*{};<2.9mm,-4mm>*{^{_2}}**@{},
   <-0.35mm,-0.35mm>*{};<-2.8mm,-4mm>*{^{_1}}**@{},
\end{xy}\Ea\hspace{-1mm}\right)
=:\hspace{-1mm}
\Ba{c}\resizebox{8mm}{!}{  \xy
(-5,6)*{}="1";
    (-5,+1)*{\circledcirc}="L";
  (-8,-5)*+{_1}="C";
   (-2,-5)*+{_2}="D";
\ar @{-} "D";"L" <0pt>
\ar @{-} "C";"L" <0pt>
\ar @{-} "1";"L" <0pt>
 \endxy}
 \Ea
$$
 which is represented
 as a $(1,2)$-corolla with the vertex denoted by $\circledcirc$; it is a cycle in the complex $(\cP,\p)$.
 Following \cite{W} one equips the extended properad $\tw \cP$ with a twisted differential $\delta$ which is defined on the generators coming from $\cP$ by
 \Beq\label{2: d_centerdot on twP under Lie}
\delta \Ba{c}\resizebox{16mm}{!}{
 \begin{xy}
 <0mm,0mm>*{\circ};<-8mm,6mm>*{^1}**@{-},
 <0mm,0mm>*{\circ};<-4.5mm,6mm>*{^2}**@{-},
 <0mm,0mm>*{\circ};<0mm,5.5mm>*{\ldots}**@{},
 <0mm,0mm>*{\circ};<3.5mm,5mm>*{}**@{-},
 <0mm,0mm>*{\circ};<8mm,6mm>*{^m}**@{-},
 <0mm,0mm>*{\circ};<-8mm,-6mm>*{_1}**@{-},
 <0mm,0mm>*{\circ};<-4.5mm,-6mm>*{_2}**@{-},
 <0mm,0mm>*{\circ};<0mm,-5.5mm>*{\ldots}**@{},
 <0mm,0mm>*{\circ};<4.5mm,-6mm>*+{}**@{-},
 <0mm,0mm>*{\circ};<8mm,-6mm>*{_n}**@{-},
   \end{xy}}\Ea
=
\p \Ba{c}\resizebox{16mm}{!}{
 \begin{xy}
 <0mm,0mm>*{\circ};<-8mm,6mm>*{^1}**@{-},
 <0mm,0mm>*{\circ};<-4.5mm,6mm>*{^2}**@{-},
 <0mm,0mm>*{\circ};<0mm,5.5mm>*{\ldots}**@{},
 <0mm,0mm>*{\circ};<3.5mm,5mm>*{}**@{-},
 <0mm,0mm>*{\circ};<8mm,6mm>*{^m}**@{-},
 <0mm,0mm>*{\circ};<-8mm,-6mm>*{_1}**@{-},
 <0mm,0mm>*{\circ};<-4.5mm,-6mm>*{_2}**@{-},
 <0mm,0mm>*{\circ};<0mm,-5.5mm>*{\ldots}**@{},
 <0mm,0mm>*{\circ};<4.5mm,-6mm>*+{}**@{-},
 <0mm,0mm>*{\circ};<8mm,-6mm>*{_n}**@{-},
   \end{xy}}\Ea
+
\overset{m-1}{\underset{i=0}{\sum}}
\Ba{c}\resizebox{17mm}{!}{
\begin{xy}
 %<0mm,0mm>*{\circ};<0mm,0mm>*{}**@{},
 <0mm,0mm>*{\circ};<-8mm,5mm>*{}**@{-},
 <0mm,0mm>*{\circ};<-3.5mm,5mm>*{}**@{-},
 <0mm,0mm>*{\circ};<-6mm,5mm>*{..}**@{},
 <0mm,0mm>*{\circ};<0mm,5mm>*{}**@{-},
  <0mm,13mm>*{\circledcirc};
  <0mm,12mm>*{};<5mm,9mm>*{\bullet}**@{-},
  <0mm,5mm>*{};<0mm,12mm>*{}**@{-},
  <0mm,14mm>*{};<0mm,17mm>*{}**@{-},
  <0mm,5mm>*{};<0mm,19mm>*{^{i\hspace{-0.2mm}+\hspace{-0.5mm}1}}**@{},
<0mm,0mm>*{\circ};<8mm,5mm>*{}**@{-},
<0mm,0mm>*{\circ};<3.5mm,5mm>*{}**@{-},
<6mm,5mm>*{..}**@{},
<-8.5mm,5.5mm>*{^1}**@{},
<-4mm,5.5mm>*{^i}**@{},
<9.0mm,5.5mm>*{^m}**@{},
 <0mm,0mm>*{\circ};<-8mm,-5mm>*{}**@{-},
 <0mm,0mm>*{\circ};<-4.5mm,-5mm>*{}**@{-},
 <-1mm,-5mm>*{\ldots}**@{},
 <0mm,0mm>*{\circ};<4.5mm,-5mm>*{}**@{-},
 <0mm,0mm>*{\circ};<8mm,-5mm>*{}**@{-},
<-8.5mm,-6.9mm>*{^1}**@{},
<-5mm,-6.9mm>*{^2}**@{},
<4.5mm,-6.9mm>*{^{n\hspace{-0.5mm}-\hspace{-0.5mm}1}}**@{},
<9.0mm,-6.9mm>*{^n}**@{},
 \end{xy}}\Ea
 - (-1)^{|a|}
\overset{n-1}{\underset{i=0}{\sum}}
 \Ba{c}\resizebox{17mm}{!}{\begin{xy}
 %<0mm,0mm>*{\circ};
 <0mm,0mm>*{\circ};<-8mm,-5mm>*{}**@{-},
 <0mm,0mm>*{\circ};<-3.5mm,-5mm>*{}**@{-},
 <0mm,0mm>*{\circ};<-6mm,-5mm>*{..}**@{},
 <0mm,0mm>*{\circ};<0mm,-5mm>*{}**@{-},
   <0mm,-11mm>*{\circledcirc};
  <0mm,-12mm>*{};<5mm,-16mm>*{\bullet}**@{-},
  <0mm,-5mm>*{};<0mm,-10mm>*{}**@{-},
  <0mm,-12mm>*{};<0mm,-17mm>*{}**@{-},
  <0mm,-5mm>*{};<0mm,-19mm>*{^{i\hspace{-0.2mm}+\hspace{-0.5mm}1}}**@{},
<0mm,0mm>*{\circ};<8mm,-5mm>*{}**@{-},
<0mm,0mm>*{\circ};<3.5mm,-5mm>*{}**@{-},
 <6mm,-5mm>*{..}**@{},
<-8.5mm,-6.9mm>*{^1}**@{},
<-4mm,-6.9mm>*{^i}**@{},
<9.0mm,-6.9mm>*{^n}**@{},
 <0mm,0mm>*{\circ};<-8mm,5mm>*{}**@{-},
 <0mm,0mm>*{\circ};<-4.5mm,5mm>*{}**@{-},
<-1mm,5mm>*{\ldots}**@{},
 <0mm,0mm>*{\circ};<4.5mm,5mm>*{}**@{-},
 <0mm,0mm>*{\circ};<8mm,5mm>*{}**@{-},
<-8.5mm,5.5mm>*{^1}**@{},
<-5mm,5.5mm>*{^2}**@{},
<4.5mm,5.5mm>*{^{m\hspace{-0.5mm}-\hspace{-0.5mm}1}}**@{},
<9.0mm,5.5mm>*{^m}**@{},
 \end{xy}}\Ea,
\Eeq
 and on the extra generator by
 \Beq\label{2: d_c on boxdot}
 \delta
 \Ba{c}\resizebox{1.7mm}{!}{\begin{xy}
 <0mm,0.5mm>*{};<0mm,4.5mm>*{}**@{-},
 <0mm,0mm>*{\bullet};<0mm,0mm>*{}**@{},
 \end{xy}}\Ea
 =
\frac{1}{2}
\Ba{c}\resizebox{8mm}{!}{  \xy
(-5,6)*{}="1";
    (-5,+1)*{\circledcirc}="L";
  (-8,-4)*{\bullet}="C";
   (-2,-4)*{\bullet}="D";
\ar @{-} "D";"L" <0pt>
\ar @{-} "C";"L" <0pt>
\ar @{-} "1";"L" <0pt>
%
 %
 %\ar @{-} "l1";"L" <0pt>
 \endxy}
 \Ea.
 \Eeq
This construction gives us the {\it twisting}\, endofunctor
$$
\tw: (\cP, \p) \lon (\tw\cP, \delta)
$$
in the category of properads under $\Lie_d$ which has been introduced first by Thomas Willwacher
in \cite{W}. This twisting endofuctor has many nice properties and applications which have been studied in \cite{W,DW,DSV}.

\sip

The distinguished cycle  $\Ba{c}\resizebox{8mm}{!}{  \xy
(-5,6)*{}="1";
    (-5,+1)*{\circledcirc}="L";
  (-8,-5)*+{_1}="C";
   (-2,-5)*+{_2}="D";
\ar @{-} "D";"L" <0pt>
\ar @{-} "C";"L" <0pt>
\ar @{-} "1";"L" <0pt>
 \endxy}
 \Ea$ in $(\cP, \p)$ remains a cycle as an element of $(\tw\cP, \delta)$ so that
the original morphism (\ref{2: i from Lie to P}) extends to the twisted version,
$$
\imath: \Lie_d \lon \tw\cP,
$$
and is given on the generator of $\Lie_d$ by the same formula.

\subsubsection{\bf From MC elements to representations of $\tw\cP$ \cite{W}}\label{2: subsec on repr of twP}
Assume we have a representation
$$
\rho: \cP \rar \cE nd_V,
$$
of  a  dg properad $\cP$ under $\Lie_d$. Composing $\rho$ with the map (\ref{2: i from Lie to P}) one makes the dg vector space $(V,d)$ into a $\Lie_d$-algebra so that it makes sense to talk about the Maurer-Cartan elements (\ref{2: MC eqn}) in $V$. Given such an element $\ga\in V$, there is
\Bi
\item[(i)] an associated twisted differential $d_\ga=d + [\ga,\ ]$ in $V$,
\item[(ii)]
 an associated representation of the twisted properad $\tw\cP$ in the dg vector space $(V,d_\ga)$,
$$
\rho^\ga: \tw\cP \lon \cE nd_V
$$
 which coincides with $\rho$ on the generators coming from $\cP$ while on the extra generator one has
$$
\rho^\ga\left(\Ba{c}\resizebox{1.6mm}{!}{\begin{xy}
 <0mm,0.5mm>*{};<0mm,4mm>*{}**@{-},
 <0mm,0mm>*{\bu};<0mm,0mm>*{}**@{},
 \end{xy}}\Ea\right):= \ga\in \cE nd_V(1,0)\equiv V.
$$
\Ei
This construction motivates the main idea of the twisting endofunctor and the terminology.

\subsection{Properad of ribbon graphs and involutive Lie bialgebras}\label{2: subsec on CycW} Let $\RGra_d=\{\RGra_d(m,n)_{m,n\geq 1}\}$ be the properad of ribbon graphs introduced
in \cite{MW}. The $\bS_m^{op}\times \bS_n$-module $\RGra_d(m,n)$ is generated by ribbon graphs $\Ga$
with $n$ labelled vertices and $m$ labelled boundaries\footnote{For a ribbon graph $\Ga$ we denote by $V(\Ga)$ its set of vertices, $B(\Ga)$ its set of boundaries and by $E(\Ga)$ its set of edges. The genus of $\Ga$ is defined by  $g= 1+\frac{1}{2}\left(\# E(\Ga) - \# V(\Ga)- \# B(\Ga)\right)$.} as, for example, the following ones,
$$
 \Ba{c}\resizebox{9mm}{!}{ \xy
 (-4,0)*{^{\bar{_1}}};
 (-1,0)*{^{\bar{_2}}};
 (1.5,0)*{^{\bar{_3}}};
  (0,5)*+{_1}*\frm{o}="1";
(0,-4)*+{_2}*\frm{o}="3";
"1";"3" **\crv{(4,0) & (4,1)};
"1";"3" **\crv{(-4,0) & (-4,-1)};
\ar @{-} "1";"3" <0pt>
\endxy}\Ea \hspace{-2mm}\in \RGra_d(3,2)
\Ba{c}\resizebox{9mm}{!}{ \xy
 (-3,1)*{^{\bar{_1}}};
 (0,8)*+{_1}*\frm{o}="1";
%(0,5)*{\circ}="1";
(0,-4)*+{_2}*\frm{o}="3";
"1";"3" **\crv{(-5,2) & (5,2)};
"1";"3" **\crv{(5,2) & (-5,2)};
"1";"3" **\crv{(-7,7) & (-7,-7)};
%\ar @{-} "1";"3" <0pt>
\endxy}\Ea \hspace{-2mm}\in \RGra_d(1,2),
 \
 \Ba{c}\resizebox{8.2mm}{!}{  \xy
 (3.5,4)*{^{\bar{1}}};
 (7,0)*+{_2}*\frm{o}="A";
 (0,0)*+{_1}*\frm{o}="B";
 \ar @{-} "A";"B" <0pt>
\endxy} \Ea\hspace{-2mm} \in \RGra(1,2),
\
 \Ba{c}\resizebox{7mm}{!}{ \xy
 (0.5,1)*{^{{^{\bar{1}}}}},
(0.5,5)*{^{{^{\bar{2}}}}},
 (0,-2)*+{_{_1}}*\frm{o}="A";
"A"; "A" **\crv{(7,7) & (-7,7)};
\endxy}\Ea \hspace{-2mm} \in \RGra_d(2,1).
$$
Such graphs are assigned the cohomological degree
$$
|\Ga|:= (1-d)\# E(\Ga)
$$
and are equipped with an {\em orientation}\, which is defined by
\Bi
\item[(i)] the ordering of edges of $\Ga$ (up to the sign action of $\bS_{\# E(\Ga)}$) for $d$ even;
\item[(ii)] the choice of a direction on each edge (up to the sign action of $\bS_2$) for $d$ odd.
\Ei
The properadic compositions are given by substituting a vertex $v$ of one ribbon graph into a boundary $b$ of another one and redistributing the edges attached to $v$ (if any) among the vertices belonging to $b$ in all possible ways while respecting the cyclic orderings (see \S 4 in \cite{MW} for full details).
\sip

A remarkable property of $\RGra_d$ is that it comes equipped with a canonical morphism of properads
\Beq\label{2: i from LBd to RGra}
i: \LBd \lon \RGra_d
\Eeq
where $\LBd$ is the properad of (degree shifted) {\it Lie bialgebras}\, defined as the quotient
of the free properad $\cF ree(\sE)$,
$$
\LB_{d}:=\cF ree\langle \sE\rangle/\langle\sR\rangle,
$$
generated by an  $\bS$-bimodule
$$
\sE(m,n):=\left\{\Ba{ll}
\id_1\ot (\sgn_2)^{\ot |d|}[d-1]=\mbox{span}\left\langle
\Ba{c}\begin{xy}
 <0mm,-0.55mm>*{};<0mm,-2.5mm>*{}**@{-},
 <0.5mm,0.5mm>*{};<2.2mm,2.2mm>*{}**@{-},
 <-0.48mm,0.48mm>*{};<-2.2mm,2.2mm>*{}**@{-},
 <0mm,0mm>*{\bu};<0mm,0mm>*{}**@{},
 %<0mm,-0.55mm>*{};<0mm,-3.8mm>*{_1}**@{},
 <0.5mm,0.5mm>*{};<2.7mm,2.8mm>*{^{_2}}**@{},
 <-0.48mm,0.48mm>*{};<-2.7mm,2.8mm>*{^{_1}}**@{},
 \end{xy}\Ea
=(-1)^{d}
\Ba{c}\begin{xy}
 <0mm,-0.55mm>*{};<0mm,-2.5mm>*{}**@{-},
 <0.5mm,0.5mm>*{};<2.2mm,2.2mm>*{}**@{-},
 <-0.48mm,0.48mm>*{};<-2.2mm,2.2mm>*{}**@{-},
 <0mm,0mm>*{\bu};<0mm,0mm>*{}**@{},
 %<0mm,-0.55mm>*{};<0mm,-3.8mm>*{_1}**@{},
 <0.5mm,0.5mm>*{};<2.7mm,2.8mm>*{^{_1}}**@{},
 <-0.48mm,0.48mm>*{};<-2.7mm,2.8mm>*{^{_2}}**@{},
 \end{xy}\Ea
   \right\rangle & \text{if}\ m=2,n=1, \vspace{3mm}\\
(\sgn_2)^{\ot |d|}\ot \id_1[d-1]=\mbox{span}\left\langle
\Ba{c}\begin{xy}
 <0mm,0.66mm>*{};<0mm,3mm>*{}**@{-},
 <0.39mm,-0.39mm>*{};<2.2mm,-2.2mm>*{}**@{-},
 <-0.35mm,-0.35mm>*{};<-2.2mm,-2.2mm>*{}**@{-},
 <0mm,0mm>*{\bu};<0mm,0mm>*{}**@{},
   <0.39mm,-0.39mm>*{};<2.9mm,-4mm>*{^{_2}}**@{},
   <-0.35mm,-0.35mm>*{};<-2.8mm,-4mm>*{^{_1}}**@{},
\end{xy}\Ea
=(-1)^{d}
\Ba{c}\begin{xy}
 <0mm,0.66mm>*{};<0mm,3mm>*{}**@{-},
 <0.39mm,-0.39mm>*{};<2.2mm,-2.2mm>*{}**@{-},
 <-0.35mm,-0.35mm>*{};<-2.2mm,-2.2mm>*{}**@{-},
 <0mm,0mm>*{\bu};<0mm,0mm>*{}**@{},
   <0.39mm,-0.39mm>*{};<2.9mm,-4mm>*{^{_1}}**@{},
   <-0.35mm,-0.35mm>*{};<-2.8mm,-4mm>*{^{_2}}**@{},
\end{xy}\Ea
\right\rangle & \text{if}\ m=1,n=2,\\
0 & \text{otherwise},
\Ea\right.
$$
by the ideal generated by the following relations
\Beq\label{2: Relations for LBd}
\sR:\left\{
\Ba{c}
\displaystyle
\oint_{123} \Ba{c}\resizebox{8.4mm}{!}{
\begin{xy}
 <0mm,0mm>*{\bu};<0mm,0mm>*{}**@{},
 <0mm,-0.49mm>*{};<0mm,-3.0mm>*{}**@{-},
 <0.49mm,0.49mm>*{};<1.9mm,1.9mm>*{}**@{-},
 <-0.5mm,0.5mm>*{};<-1.9mm,1.9mm>*{}**@{-},
 <-2.3mm,2.3mm>*{\bu};<-2.3mm,2.3mm>*{}**@{},
 <-1.8mm,2.8mm>*{};<0mm,4.9mm>*{}**@{-},
 <-2.8mm,2.9mm>*{};<-4.6mm,4.9mm>*{}**@{-},
   <0.49mm,0.49mm>*{};<2.7mm,2.3mm>*{^3}**@{},
   <-1.8mm,2.8mm>*{};<0.4mm,5.3mm>*{^2}**@{},
   <-2.8mm,2.9mm>*{};<-5.1mm,5.3mm>*{^1}**@{},
 \end{xy}}\Ea =0
 \ , \ \ \ \ \
%%%%%%%%%%%%%% Lie %%%%%%%%%%%%%%%%%%%%%%%%
\oint_{123} \Ba{c}\resizebox{8.4mm}{!}{ \begin{xy}
 <0mm,0mm>*{\bu};<0mm,0mm>*{}**@{},
 <0mm,0.69mm>*{};<0mm,3.0mm>*{}**@{-},
 <0.39mm,-0.39mm>*{};<2.4mm,-2.4mm>*{}**@{-},
 <-0.35mm,-0.35mm>*{};<-1.9mm,-1.9mm>*{}**@{-},
 <-2.4mm,-2.4mm>*{\bu};<-2.4mm,-2.4mm>*{}**@{},
 <-2.0mm,-2.8mm>*{};<0mm,-4.9mm>*{}**@{-},
 <-2.8mm,-2.9mm>*{};<-4.7mm,-4.9mm>*{}**@{-},
    <0.39mm,-0.39mm>*{};<3.3mm,-4.0mm>*{^3}**@{},
    <-2.0mm,-2.8mm>*{};<0.5mm,-6.7mm>*{^2}**@{},
    <-2.8mm,-2.9mm>*{};<-5.2mm,-6.7mm>*{^1}**@{},
 \end{xy}}\Ea =0
\vspace{2mm} \\
%%%%%%%%%%%%%%%%%%%%%%% Lie[1]Bi %%%%%%%%%%%%%%%
 \Ba{c}\resizebox{5.2mm}{!}{\begin{xy}
 <0mm,2.47mm>*{};<0mm,0.12mm>*{}**@{-},
 <0.5mm,3.5mm>*{};<2.2mm,5.2mm>*{}**@{-},
 <-0.48mm,3.48mm>*{};<-2.2mm,5.2mm>*{}**@{-},
 <0mm,3mm>*{\bu};<0mm,3mm>*{}**@{},
  <0mm,-0.8mm>*{\bu};<0mm,-0.8mm>*{}**@{},
<-0.39mm,-1.2mm>*{};<-2.2mm,-3.5mm>*{}**@{-},
 <0.39mm,-1.2mm>*{};<2.2mm,-3.5mm>*{}**@{-},
     <0.5mm,3.5mm>*{};<2.8mm,5.7mm>*{^2}**@{},
     <-0.48mm,3.48mm>*{};<-2.8mm,5.7mm>*{^1}**@{},
   <0mm,-0.8mm>*{};<-2.7mm,-5.2mm>*{^1}**@{},
   <0mm,-0.8mm>*{};<2.7mm,-5.2mm>*{^2}**@{},
\end{xy}}\Ea
+(-1)^{d}
\Ba{c}\resizebox{7mm}{!}{\begin{xy}
 <0mm,-1.3mm>*{};<0mm,-3.5mm>*{}**@{-},
 <0.38mm,-0.2mm>*{};<2.0mm,2.0mm>*{}**@{-},
 <-0.38mm,-0.2mm>*{};<-2.2mm,2.2mm>*{}**@{-},
<0mm,-0.8mm>*{\bu};<0mm,0.8mm>*{}**@{},
 <2.4mm,2.4mm>*{\bu};<2.4mm,2.4mm>*{}**@{},
 <2.77mm,2.0mm>*{};<4.4mm,-0.8mm>*{}**@{-},
 <2.4mm,3mm>*{};<2.4mm,5.2mm>*{}**@{-},
     <0mm,-1.3mm>*{};<0mm,-5.3mm>*{^1}**@{},
     <2.5mm,2.3mm>*{};<5.1mm,-2.6mm>*{^2}**@{},
    <2.4mm,2.5mm>*{};<2.4mm,5.7mm>*{^2}**@{},
    <-0.38mm,-0.2mm>*{};<-2.8mm,2.5mm>*{^1}**@{},
    \end{xy}}\Ea
  +
\Ba{c}\resizebox{7mm}{!}{\begin{xy}
 <0mm,-1.3mm>*{};<0mm,-3.5mm>*{}**@{-},
 <0.38mm,-0.2mm>*{};<2.0mm,2.0mm>*{}**@{-},
 <-0.38mm,-0.2mm>*{};<-2.2mm,2.2mm>*{}**@{-},
<0mm,-0.8mm>*{\bu};<0mm,0.8mm>*{}**@{},
 <2.4mm,2.4mm>*{\bu};<2.4mm,2.4mm>*{}**@{},
 <2.77mm,2.0mm>*{};<4.4mm,-0.8mm>*{}**@{-},
 <2.4mm,3mm>*{};<2.4mm,5.2mm>*{}**@{-},
     <0mm,-1.3mm>*{};<0mm,-5.3mm>*{^2}**@{},
     <2.5mm,2.3mm>*{};<5.1mm,-2.6mm>*{^1}**@{},
    <2.4mm,2.5mm>*{};<2.4mm,5.7mm>*{^2}**@{},
    <-0.38mm,-0.2mm>*{};<-2.8mm,2.5mm>*{^1}**@{},
    \end{xy}}\Ea
  + (-1)^{d}
\Ba{c}\resizebox{7mm}{!}{\begin{xy}
 <0mm,-1.3mm>*{};<0mm,-3.5mm>*{}**@{-},
 <0.38mm,-0.2mm>*{};<2.0mm,2.0mm>*{}**@{-},
 <-0.38mm,-0.2mm>*{};<-2.2mm,2.2mm>*{}**@{-},
<0mm,-0.8mm>*{\bu};<0mm,0.8mm>*{}**@{},
 <2.4mm,2.4mm>*{\bu};<2.4mm,2.4mm>*{}**@{},
 <2.77mm,2.0mm>*{};<4.4mm,-0.8mm>*{}**@{-},
 <2.4mm,3mm>*{};<2.4mm,5.2mm>*{}**@{-},
     <0mm,-1.3mm>*{};<0mm,-5.3mm>*{^2}**@{},
     <2.5mm,2.3mm>*{};<5.1mm,-2.6mm>*{^1}**@{},
    <2.4mm,2.5mm>*{};<2.4mm,5.7mm>*{^1}**@{},
    <-0.38mm,-0.2mm>*{};<-2.8mm,2.5mm>*{^2}**@{},
    \end{xy}}\Ea
 +
\Ba{c}\resizebox{7mm}{!}{\begin{xy}
 <0mm,-1.3mm>*{};<0mm,-3.5mm>*{}**@{-},
 <0.38mm,-0.2mm>*{};<2.0mm,2.0mm>*{}**@{-},
 <-0.38mm,-0.2mm>*{};<-2.2mm,2.2mm>*{}**@{-},
<0mm,-0.8mm>*{\bu};<0mm,0.8mm>*{}**@{},
 <2.4mm,2.4mm>*{\bu};<2.4mm,2.4mm>*{}**@{},
 <2.77mm,2.0mm>*{};<4.4mm,-0.8mm>*{}**@{-},
 <2.4mm,3mm>*{};<2.4mm,5.2mm>*{}**@{-},
     <0mm,-1.3mm>*{};<0mm,-5.3mm>*{^1}**@{},
     <2.5mm,2.3mm>*{};<5.1mm,-2.6mm>*{^2}**@{},
    <2.4mm,2.5mm>*{};<2.4mm,5.7mm>*{^1}**@{},
    <-0.38mm,-0.2mm>*{};<-2.8mm,2.5mm>*{^2}**@{},
    \end{xy}}\Ea=0.
    \Ea
\right.
\Eeq
The case $d=1$ corresponds to the standard notion of Lie bialgebra introduced by V.\ Drinfeld in \cite{Dr1}.

\sip

The morphism (\ref{2: i from LBd to RGra}) is given on the generators by \cite{MW}
\Beq\label{2: values on generators  of i from LBd to RGra}
i\left(
\Ba{c}\begin{xy}
 <0mm,0.66mm>*{};<0mm,3mm>*{}**@{-},
 <0.39mm,-0.39mm>*{};<2.2mm,-2.2mm>*{}**@{-},
 <-0.35mm,-0.35mm>*{};<-2.2mm,-2.2mm>*{}**@{-},
 <0mm,0mm>*{\bu};<0mm,0mm>*{}**@{},
   <0.39mm,-0.39mm>*{};<2.9mm,-4mm>*{^{_2}}**@{},
   <-0.35mm,-0.35mm>*{};<-2.8mm,-4mm>*{^{_1}}**@{},
 <0mm,4mm>*{^{_{\bar{1}}}}**@{},
\end{xy}\Ea
\right):=\Ba{c}\resizebox{8.2mm}{!}{  \xy
 (3.5,4)*{^{\bar{1}}};
 (7,0)*+{_2}*\frm{o}="A";
 (0,0)*+{_1}*\frm{o}="B";
 \ar @{-} "A";"B" <0pt>
\endxy} \Ea,
\ \ \ \ \
i\left(\Ba{c}\begin{xy}
 <0mm,-0.55mm>*{};<0mm,-2.5mm>*{}**@{-},
 <0.5mm,0.5mm>*{};<2.2mm,2.2mm>*{}**@{-},
 <-0.48mm,0.48mm>*{};<-2.2mm,2.2mm>*{}**@{-},
 <0mm,0mm>*{\bu};<0mm,0mm>*{}**@{},
 <0mm,-0.55mm>*{};<0mm,-3.8mm>*{_1}**@{},
 <0.5mm,0.5mm>*{};<2.7mm,2.8mm>*{^{_{\bar{1}}}}**@{},
 <-0.48mm,0.48mm>*{};<-2.7mm,2.8mm>*{^{_{\bar{2}}}}**@{},
 \end{xy}\Ea\right)=  \Ba{c}\resizebox{7mm}{!}{ \xy
 (0.5,1)*{^{{^{\bar{1}}}}},
(0.5,5)*{^{{^{\bar{2}}}}},
 (0,-2)*+{_{_1}}*\frm{o}="A";
"A"; "A" **\crv{(7,7) & (-7,7)};
\endxy}\Ea
\Eeq

The properad $\LoBd$ of {\it involutive}\, Lie bialgebras is the quotient of $\LBd$ by the ideal generated by the following relation
$$
\Ba{c}\resizebox{5mm}{!}
{\xy
 (0,0)*{\bu}="a",
(0,6)*{\bu}="b",
(3,3)*{}="c",
(-3,3)*{}="d",
 (0,9)*{}="b'",
(0,-3)*{}="a'",
%"a";"c"**\dir{-};
\ar@{-} "a";"c" <0pt>
\ar @{-} "a";"d" <0pt>
\ar @{-} "a";"a'" <0pt>
\ar @{-} "b";"c" <0pt>
\ar @{-} "b";"d" <0pt>
\ar @{-} "b";"b'" <0pt>
\endxy}
\Ea=0.
$$
The map (\ref{2: i from LBd to RGra}) factors through the canonical projection
$$
\LBd \lon \LoBd
$$
so that the formulae (\ref{2: values on generators  of i from LBd to RGra}) define in fact a
canonical morphism (denoted by the same symbol $i$) \cite{MW}
$$
i: \LoBd \lon \RGra_d.
$$
The morphism $i$ was used in \cite{MW} to built several ribbon graph complexes with interesting cohomology; in this paper we consider its applications in string topology.
 The minimal resolution of $\LoBd$ has been constructed in \cite{CMW} and is denoted by $\HoLoBd$; this is a free properad generated
by the following family of (skew)symmetric corollas of degrees $1-d(m+n+2a-2)$,
\Beq\label{2: generators of HoLoB}
\left\{
\Ba{c}\resizebox{16mm}{!}{\xy
(-9,-6)*{};
(0,0)*+{a}*\cir{}
**\dir{-};
(-5,-6)*{};
(0,0)*+{a}*\cir{}
**\dir{-};
(9,-6)*{};
(0,0)*+{a}*\cir{}
**\dir{-};
(5,-6)*{};
(0,0)*+{a}*\cir{}
**\dir{-};
(0,-6)*{\ldots};
(-10,-8)*{_1};
(-6,-8)*{_2};
(10,-8)*{_n};
(-9,6)*{};
(0,0)*+{a}*\cir{}
**\dir{-};
(-5,6)*{};
(0,0)*+{a}*\cir{}
**\dir{-};
(9,6)*{};
(0,0)*+{a}*\cir{}
**\dir{-};
(5,6)*{};
(0,0)*+{a}*\cir{}
**\dir{-};
(0,6)*{\ldots};
(-10,8)*{_1};
(-6,8)*{_2};
(10,8)*{_m};
\endxy}\Ea
=(-1)^{d(\sigma+\tau)}
\Ba{c}\resizebox{20mm}{!}{\xy
(-9,-6)*{};
(0,0)*+{a}*\cir{}
**\dir{-};
(-5,-6)*{};
(0,0)*+{a}*\cir{}
**\dir{-};
(9,-6)*{};
(0,0)*+{a}*\cir{}
**\dir{-};
(5,-6)*{};
(0,0)*+{a}*\cir{}
**\dir{-};
(0,-6)*{\ldots};
(-12,-8)*{_{\tau(1)}};
(-6,-8)*{_{\tau(2)}};
(12,-8)*{_{\tau(n)}};
(-9,6)*{};
(0,0)*+{a}*\cir{}
**\dir{-};
(-5,6)*{};
(0,0)*+{a}*\cir{}
**\dir{-};
(9,6)*{};
(0,0)*+{a}*\cir{}
**\dir{-};
(5,6)*{};
(0,0)*+{a}*\cir{}
**\dir{-};
(0,6)*{\ldots};
(-12,8)*{_{\sigma(1)}};
(-6,8)*{_{\sigma(2)}};
(12,8)*{_{\sigma(m)}};
\endxy}\Ea\ \ \ \forall \sigma\in \bS_m, \forall \tau\in \bS_n\right\}_{m+n+ a\geq 3 \atop m\geq 1, n\geq 1, a\geq 0}
\Eeq
The differential
 is given by
\Beq\label{2: d on Lie inv infty}
\delta
\Ba{c}\resizebox{16mm}{!}{\xy
(-9,-6)*{};
(0,0)*+{a}*\cir{}
**\dir{-};
(-5,-6)*{};
(0,0)*+{a}*\cir{}
**\dir{-};
(9,-6)*{};
(0,0)*+{a}*\cir{}
**\dir{-};
(5,-6)*{};
(0,0)*+{a}*\cir{}
**\dir{-};
(0,-6)*{\ldots};
(-10,-8)*{_1};
(-6,-8)*{_2};
(10,-8)*{_n};
(-9,6)*{};
(0,0)*+{a}*\cir{}
**\dir{-};
(-5,6)*{};
(0,0)*+{a}*\cir{}
**\dir{-};
(9,6)*{};
(0,0)*+{a}*\cir{}
**\dir{-};
(5,6)*{};
(0,0)*+{a}*\cir{}
**\dir{-};
(0,6)*{\ldots};
(-10,8)*{_1};
(-6,8)*{_2};
(10,8)*{_m};
\endxy}\Ea
=
\sum_{l\geq 1}\sum_{a=b+c+l-1}\sum_{[m]=I_1\sqcup I_2\atop
[n]=J_1\sqcup J_2} \pm
\Ba{c}
%
%
%%%%%%%%%%%%%%%% two vertex graph with l internal edges %%%%%%%%%%
\Ba{c}\resizebox{21mm}{!}{\xy
(0,0)*+{b}*\cir{}="b",
(10,10)*+{c}*\cir{}="c",
%
%%%%%%%%%% edges to b %%%%%%%%%%%%
(-9,6)*{}="1",
(-7,6)*{}="2",
(-2,6)*{}="3",
(-3.5,5)*{...},
(-4,-6)*{}="-1",
(-2,-6)*{}="-2",
(4,-6)*{}="-3",
(1,-5)*{...},
(0,-8)*{\underbrace{\ \ \ \ \ \ \ \ }},
(0,-11)*{_{J_1}},
(-6,8)*{\overbrace{ \ \ \ \ \ \ }},
(-6,11)*{_{I_1}},
%%%%%%%%%% edges to c %%%%%%%%%%%%
(6,16)*{}="1'",
(8,16)*{}="2'",
(14,16)*{}="3'",
(11,15)*{...},
(11,6)*{}="-1'",
(16,6)*{}="-2'",
(18,6)*{}="-3'",
(13.5,6)*{...},
(15,4)*{\underbrace{\ \ \ \ \ \ \ }},
(15,1)*{_{J_2}},
(10,18)*{\overbrace{ \ \ \ \ \ \ \ \ }},
(10,21)*{_{I_2}},
%
%%%%%%%%%%% internal curved edges %%%%%%%%%%%
(0,2)*-{};(8.0,10.0)*-{}
**\crv{(0,10)};
(0.5,1.8)*-{};(8.5,9.0)*-{}
**\crv{(0.4,7)};
(1.5,0.5)*-{};(9.1,8.5)*-{}
**\crv{(5,1)};
(1.7,0.0)*-{};(9.5,8.6)*-{}
**\crv{(6,-1)};
(5,5)*+{...};
\ar @{-} "b";"1" <0pt>
\ar @{-} "b";"2" <0pt>
\ar @{-} "b";"3" <0pt>
\ar @{-} "b";"-1" <0pt>
\ar @{-} "b";"-2" <0pt>
\ar @{-} "b";"-3" <0pt>
\ar @{-} "c";"1'" <0pt>
\ar @{-} "c";"2'" <0pt>
\ar @{-} "c";"3'" <0pt>
\ar @{-} "c";"-1'" <0pt>
\ar @{-} "c";"-2'" <0pt>
\ar @{-} "c";"-3'" <0pt>
\endxy}\Ea
%%%%%%%%%%%%%%%%%%%
\Ea
\Eeq
where the summation parameter $l$ counts the number of internal edges connecting the two vertices
on the r.h.s.

%%%%%%%%%%%%%%%%%%%%%%%%%%%%%%%%%%%%%%%%%%%%%%%%%%%%%%%%%%%%%%%%%%%%%%%%%%%%%%%%

\subsection{Gravity properad} The morphism (\ref{2: i from LBd to RGra}) says, in particular, that $\RGra_d$ is a properad under $\Lie_d$ and hence can be twisted as explained in \S 2.1. The twisted properad  $\tw\cR\cG ra_d=\{\tw\cR\cG ra_d(m,n)\}_{m\geq 1, n\geq 0}$ is generated by ribbon graphs $\Ga$ with all boundaries labelled and
with vertices of two types --- the white vertices which are labelled, and the black vertices which are unlabelled but assigned the cohomological degree $d$; their orientation is defined as in the case of ribbon graphs from $\RGra_d$ except that for $d$ odd one has to choose in addition an ordering (defined up to the sign action of the permutation group) of black vertices. A ribbon graph $\Ga\in \tw\RGra_d(m,n)$ has $m$ labelled boundaries, $n$ labelled white vertices and any finite number $k\geq 0$ of black vertices; its cohomological degree is defined by
$$
|\Ga|=(1-d)\# E(\Ga) + kd
$$
 For example,
$$
\stackrel{^{\bar{1}}}{\bu}\ \in \tw\RGra_d(1,0),\ \ \
\Ba{c}\resizebox{9mm}{!}{ \xy
 (-3,1)*{^{\bar{_1}}};
 (0,8)*+{_1}*\frm{o}="1";
%(0,5)*{\circ}="1";
(0,-4)*{\bu}="3";
"1";"3" **\crv{(-5,2) & (5,2)};
"1";"3" **\crv{(5,2) & (-5,2)};
"1";"3" **\crv{(-7,7) & (-7,-7)};
%\ar @{-} "1";"3" <0pt>
\endxy}\Ea \in \tw\cR\cG ra_d(1,1),
 \ \ \
 \Ba{c}\resizebox{7mm}{!}{  \xy
 (3.1,4)*{^{\bar{1}}};
 (6,0)*{\bu}="A";
 (0,0)*{\bu}="B";
 \ar @{-} "A";"B" <0pt>
\endxy} \Ea \in \tw\cR\cG ra_d(1,0),
\ \
 \xy
 (0.5,1)*{^{{^{\bar{1}}}}},
(0.5,5)*{^{{^{\bar{2}}}}},
 (0,-2)*+{_{_1}}*\frm{o}="A";
% (0,0)*+{_1}*\frm{o}="B";
%(0,-2)*{\circ}="A";
%(0,-2)*{\circ}="B";
"A"; "A" **\crv{(7,7) & (-7,7)};
\endxy\in \tw\cR\cG ra_d(2,1),
%\xy
%
% (0,0)*{\bu}="a",
%(6,0)*{\bu}="b",
%
%\ar @{-} "a";"b" <0pt>
%\endxy \in \Tw\cR\cG ra_d(1,0).
$$
The first graph in this list is precisely the ribbon graph incarnation of the extra $(1,0)$-generator
which --- in the decorated corolla notation --- is represented  in \S 2.1 as $\Ba{c}\resizebox{1.7mm}{!}{\begin{xy}
 <0mm,0.5mm>*{};<0mm,3.5mm>*{}**@{-},
 <0mm,0mm>*{\bullet};<0mm,0mm>*{}**@{},
 \end{xy}}\Ea$. Substituting the unique boundary of the ribbon graph $\bu$ into a white vertex of graphs from $\RGra_d$ one obtains generic elements of $\tw \RGra_d$ as described above. Given two ribbon graphs
 $\Ga_1,\Ga_2\in \tw \RGra_d$, we denote by $\Ga_1\, {}_i\hspace{-0.3mm}\circ_j \Ga_2$ their partial properadic composition
 given by substituting the $j$-th boundary of $\Ga_2$ into the $i$-th white vertex of $\Ga_1$. The set of black (resp., white) vertices of $\Ga\in \tw\RGra_d$ is denoted by $V_\bu(\Ga)$ (resp., $V_\circ(\Ga)$).

 \sip

 As the properad $\RGra_d$ has trivial differential, the induced differential $\delta$ in
  $\tw \RGra_d$ is completely determined by the last two summands in (\ref{2: d_centerdot on twP under Lie}) which in the ribbon graph incarnation give us the following explicit formula
\Beq\label{2: delta in Tw(RGra)}
\delta\Ga:=
\sum_{i=1}^m \Ba{c}\resizebox{3mm}{!}{  \xy
 (0,6)*{\bu}="A";
 (0,0)*+{_1}*\frm{o}="B";
 \ar @{-} "A";"B" <0pt>
\endxy} \Ea
\  _1\hspace{-0.7mm}\circ_j \Ga\ \
- \ \ (-1)^{|\Ga|} %\sum_{v\in V_\circ(\Ga)}
\sum_{j=1}^n\Ga\  _j\hspace{-0.7mm}\circ_1
\Ba{c}\resizebox{3mm}{!}{  \xy
 (0,6)*{\bu}="A";
 (0,0)*+{_1}*\frm{o}="B";
 \ar @{-} "A";"B" <0pt>
\endxy} \Ea
\  -(-1)^{|\Ga|}\ \frac{1}{2} \sum_{v\in V_\bu(\Ga)} \Ga\circ_v  \left(\xy
 (0,0)*{\bullet}="a",
(5,0)*{\bu}="b",
\ar @{-} "a";"b" <0pt>
\endxy\right)
\Eeq
where the symbol $\Ga\circ_v  \left(\xy
 (0,0)*{\bullet}="a",
(5,0)*{\bu}="b",
\ar @{-} "a";"b" <0pt>
\endxy\right)$ means substituting the unique boundary of the ribbon graph  $\xy
 (0,0)*{\bullet}="a",
(5,0)*{\bu}="b",
\ar @{-} "a";"b" <0pt>
\endxy$ into the black vertex $v$ of the graph $\Ga$ and then taking the sum over all possible re-attachments of the half-edges attached earlier  to $v$ among the two newly created black vertices in a way which respects the cyclic orderings. Note that  if $\Ga$ has no univalent black vertices, then $\delta \Ga$ has no univalent black vertices as well (they all cancel out when applying the above formula).

\sip

By the general rule,
the twisted properad $\tw \RGra_d$ comes equipped with a morphism
$$
\Ba{rccc}
i: & \Lie_d & \lon & \tw \RGra_d\\
   & \Ba{c}\begin{xy}
 <0mm,0.66mm>*{};<0mm,3mm>*{}**@{-},
 <0.39mm,-0.39mm>*{};<2.2mm,-2.2mm>*{}**@{-},
 <-0.35mm,-0.35mm>*{};<-2.2mm,-2.2mm>*{}**@{-},
 <0mm,0mm>*{\bu};<0mm,0mm>*{}**@{},
   <0.39mm,-0.39mm>*{};<2.9mm,-4mm>*{^{_2}}**@{},
   <-0.35mm,-0.35mm>*{};<-2.8mm,-4mm>*{^{_1}}**@{},
 <0mm,4mm>*{^{_{\bar{1}}}}**@{},
\end{xy}\Ea
&\lon&
\Ba{c}\resizebox{8.2mm}{!}{  \xy
 (3.5,4)*{^{\bar{1}}};
 (7,0)*+{_2}*\frm{o}="A";
 (0,0)*+{_1}*\frm{o}="B";
 \ar @{-} "A";"B" <0pt>
\endxy} \Ea
\Ea
$$
induced by the above morphism (\ref{2: i from LBd to RGra}) from $\Lie_d$ to $\RGra_d$. However that morphism does not extend  to a morphism from $\LBd$ to $\tw \RGra_d$ as the image of the second generator of $\LBd$  is no more a cocycle \cite{Me2},
\Beq\label{2: delta on ribbon CoLie corolla}
\delta \xy
 (0.5,1)*{^{{^{\bar{1}}}}},
(0.5,5)*{^{{^{\bar{2}}}}},
 (0,-2)*+{_{_1}}*\frm{o}="A";
% (0,0)*+{_1}*\frm{o}="B";
%(0,-2)*{\circ}="A";
%(0,-2)*{\circ}="B";
"A"; "A" **\crv{(7,7) & (-7,7)};
\endxy=
\Ba{c}\resizebox{6mm}{!}{
\mbox{$\xy
 (0.5,0.9)*{^{{^{\bar{1}}}}},
(0.5,5)*{^{{^{\bar{2}}}}},
 (0,-8)*+{_{_1}}*\frm{o}="C";
(0,-2)*{\bu}="A";
(0,-2)*{\bu}="B";
"A"; "B" **\crv{(6,6) & (-6,6)};
 \ar @{-} "A";"C" <0pt>
\endxy$}}
\Ea
+ (-1)^d
\Ba{c}\resizebox{6mm}{!}{
\mbox{$\xy
 (0.5,0.9)*{^{{^{\bar{2}}}}},
(0.5,5)*{^{{^{\bar{1}}}}},
 (0,-8)*+{_{_1}}*\frm{o}="C";
(0,-2)*{\bu}="A";
(0,-2)*{\bu}="B";
"A"; "B" **\crv{(6,6) & (-6,6)};
 \ar @{-} "A";"C" <0pt>
\endxy$}}
\Ea \neq 0
\Eeq
 Nevertheless the cohomology properad $H^\bu(\tw\RGra_d,\delta)$ is proven in \cite{Me2} to be a properad under the properad of degree shifted {\it quasi}\, Lie bialgebras; remarkably, the latter properad takes care about {\it infinitely many}\, cohomology classes in $H^\bu(\tw\RGra_d,\delta)$ (see \S 3.9 in \cite{Me2}). We do not use this fact in this paper and hence omit any details.

\sip

One of the main results in \cite{Me2} is the computation of the cohomology properad $H^\bu(\tw\RGra_d,\delta)$  using the remarkable geometric  K.\ Costello's theory \cite{Co1,Co2} of partially compactified moduli spaces $\overline{\cN}_{g,m,0,n}$ of connected stable Riemann
surfaces $\Sigma$ of genus $g$ with $m$ marked
boundary components and $n$ marked points
in the interior $\Sigma\setminus \p\Sigma$ and with possible nodes on the boundary $\p\Sigma$.
To formulate it precisely we notice that the differential $\delta$ in each complex $\tw\RGra_d(m,n)$ preserves the genus of ribbon graphs so that the latter decomposes into a direct sum,
$$
\tw\RGra_d(m,n)=\bigoplus_{g\geq 0} \tw\RGra_d(g;m,n)
$$
where the subcomplex $\tw\RGra_d(g;m,n)$ is spanned by ribbon graphs of genus $g$. Let $\cM_{g,N}$ stand for the moduli space of algebraic curves $\Sigma$ of genus $g$ with $N\geq 1$ marked points. The group $\bS_N$ acts on the compactly supported cohomology group $H^\bu_c(\cM_{g,N})$ by relabeling. Assume that the set $S\subset \Sigma$ of $N$ marked points in an algebraic curve $\Sigma$ is split into the disjoint union of subsets of (so called) ``in"- and ``out" points,
$$
S=S_{in}\sqcup S_{out}
$$
which have cardinalities $n:=\# S_{in}\geq 0$ and $m:=\# S_{out}\geq 1$, and denote the moduli space of such algebraic curves by $\cM_{g,m+n}$; its compactly supported cohomology group $H_c^\bu(\cM_{g,m+n})$ is naturally an $\bS^{op}_m\times \bS_n$-module. The tensor product
 $$
 H^{\bu}_c(\cM_{g,m+n})\ot \sgn_m
 $$
is also an $\bS^{op}_m\times \bS_n$-module with group $\bS_m^{op}$ acting, by definition, diagonally.

\subsubsection{\bf Theorem \cite{Me2}}
(i) {\em For any $g\geq 0$, $m\geq 1$ and $n\geq 0$ with $2g+m+n\geq 3$ one has an isomorphism of $\bS_m^{op}\times \bS_n$-modules,}
$$
H^\bu(\tw\RGra_d(g;m,n))= H^{\bu-m +d(2g-2+m+n)}_c(\cM_{g,m+n})\ot \sgn_m\simeq H^{\bu+d(2g-2+m+n)}_c(\cM_{g,m+n}\times \R^m)
$$
(ii) {\em For any $g\geq 0$, $m\geq 1$ and $n\geq 0$ with $2g+m+n< 3$ one has}
$$
H^k(\tw\RGra_d(g;m,n))=\left\{\Ba{ll} \K & \text{if} \ g=0,n=0,m=2, k\geq 1 \ \& \ k\equiv 2d+1 \bmod 4 \\
0 & \text{otherwise}.
\Ea
\right.
$$
{\em where $\K$ is generated by  the unique polytope-like ribbon graph with $k$ edges and $k$ bivalent vertices which are all black}.

\sip

The polytope-like cohomology classes in $H^\bu(\tw\RGra_d)$ play no role in this paper so that we can restrict to the dg sub-properad
$$
\cC h\Grav_d \subset \tw\RGra_d
$$
spanned by ribbon graphs whose black vertices are at least {\em trivalent}; it is called the {\em chain gravity properad}\, in \cite{Me2} and was denoted there simply by $\tw^{\geq 3}\RGra_d$. Its cohomology properad is called the {\em gravity properad}\, and is given by (\ref{1: def of Grav_d});
 it contains E.\ Getzler's gravity operad \cite{Ge} as the genus zero sub-operad (\ref{1: grav_d operad}).
% genus zero sub-properad
%\Beq\label{2: grav}
%grav_d=\{ H^{\bu-1+d(n-1)}_c(\cM_{0,1+n})\}_{n\geq 2};
%\Eeq
%the proof of this statement given in \cite{Me2} is based on the work of B.\ Ward %\cite{Wa}.

\bip

\bip

%%%%%%%%%%%%%%%%%%%%%%%%%%%%%%%%%%%%%%%%%%%%%%%%%%%%%%%%%%%%%%%%%%%%%%%%%%%%
%%%%%%%%%%%%%%%%%%%%%%%%%%%%%%%%%%%%%%%%%%%%%%%%%%%%%%%%%%%%%%%%%%%%%%%%%%%%

{\Large
\section{\bf From gravity to string topology}
}

\mip

\subsection{A class of representations of $\RGra_d$}  The properadic compositions in $\RGra_d$ have been designed in \cite{MW} in such a way that $\RGra_d$ admits a canonical representation
\Beq\label{3: rho from RGra to CycW}
\rho: \RGra_d \lon \cE nd_{Cyc(W)},
\Eeq
 in the space of ``cyclic words",
 $$
 Cyc(W):= \bigoplus_{k\geq 0} \left(\ot^k W\right)_{\Z^*_k},
 $$
which is generated by a graded vector space $W$  equipped with a scalar product of degree $1-d$,
$$
\Ba{rccc}
\Theta: & W\ot W & \lon & \K[1-d] \\
  &  w_1\ot w_2 & \lon & \Theta(w_1,w_2)
\Ea
$$
such that
\Beq\label{3: skewsymmetry on scalar product}
\Theta(w_2,w_1)=(-1)^{d+|w_1||w_2|}\Theta(w_1,w_2).
\Eeq
%The idea behind the above representation map $\rho$ is best described in a concrete %example. Consider a ribbon graph
%$$
%\Ga=
%$$
%\ldots

\sip

We need in this paper a slightly sharpened version of the morphism (\ref{3: rho from RGra to CycW}). Note that a differential $\p$ in $W$ induces a differential in $Cyc(W)$ and hence in the endomorphism properad $\cE nd_{Cyc(W)}$ which is denoted by the same letter $\p$.

\subsubsection{\bf Proposition}\label{3: Prop on repr of RGra_2} {\em Let $(W,\p)$ be as dg vector space equipped with a scalar product $\Theta:
W\ot W \rar \K[1-d]$ which satisfies (\ref{3: skewsymmetry on scalar product}) and which is compatible with the differential in the sense that
$$
\Theta(\p w_1, w_2) + (-1)^{|w_1|}\Theta(w_1,\p w_2)=0\ \ \ \ \forall\ w_1,w_2\in W.
$$
Then the morphism (\ref{3: rho from RGra to CycW}) extends to a morphism of {\em differential}\, properads
\Beq\label{3: rho from RGra to differential CycW}
\rho: (\RGra_d,0) \lon (\cE nd_{Cyc(W)}, \p).
\Eeq
where $\RGra_d$ is equipped with the trivial differential.}

\mip

The proposition is equivalent to saying that for any $\Ga\in \RGra_d$ its image $\rho(\Ga)$ is a cycle in the complex
$(\cE nd_{Cyc(W)}, \p)$. Which is easy to check using the compatibility of $\Theta$ with $\p$ and the explicit formula for $\rho(\Ga)$ given in \S 4.2.2 of \cite{MW}.

%(Hence
%$$
%(-1)^{d+ |w_1||w_2|+ |w_2|}\Theta(w_2,\p w_1) + (-1)^{|w_1|}(-1)^{d+ |w_1||w_2|+|w_1|}\Theta(\p w_2, %w_1)=0
%$$
%OK.)

\sip

Composing the canonical morphism
\Beq\label{3: morphism from LoB to RGra}
\Ba{rccc}
i: & \LoB_{d} & \lon & \cR\cG ra_{d}\vspace{2mm}\\
&
\Ba{c}\begin{xy}
 <0mm,-0.55mm>*{};<0mm,-2.5mm>*{}**@{-},
 <0.5mm,0.5mm>*{};<2.2mm,2.5mm>*{}**@{-},
 <-0.48mm,0.48mm>*{};<-2.2mm,2.5mm>*{}**@{-},
 <0mm,0mm>*{\bu};<0mm,0mm>*{}**@{},
 <0mm,-0.55mm>*{};<0mm,-3.8mm>*{_{_1}}**@{},
 <0.5mm,0.5mm>*{};<2.7mm,3.1mm>*{^{^{\bar{2}}}}**@{},
 <-0.48mm,0.48mm>*{};<-2.7mm,3.1mm>*{^{^{\bar{1}}}}**@{},
 \end{xy}\Ea
 &\lon &
 \Ba{c}\resizebox{7mm}{!}{ \xy
 (0.5,1)*{^{{^{\bar{1}}}}},
(0.5,5)*{^{{^{\bar{2}}}}},
 (0,-2)*+{_{_1}}*\frm{o}="A";
"A"; "A" **\crv{(7,7) & (-7,7)};
\endxy}\Ea
\vspace{2mm}\\
&
\Ba{c}\begin{xy}
 <0mm,0.66mm>*{};<0mm,3mm>*{}**@{-},
 <0.39mm,-0.39mm>*{};<2.2mm,-2.2mm>*{}**@{-},
 <-0.35mm,-0.35mm>*{};<-2.2mm,-2.2mm>*{}**@{-},
 <0mm,0mm>*{\bu};<0mm,4.1mm>*{^{^{\bar{1}}}}**@{},
   <0.39mm,-0.39mm>*{};<2.9mm,-4mm>*{^{_2}}**@{},
   <-0.35mm,-0.35mm>*{};<-2.8mm,-4mm>*{^{_1}}**@{},
\end{xy}\Ea
 & \lon &
 \Ba{c}\resizebox{10mm}{!}{  \xy
 (3.5,4)*{^{\bar{1}}};
 (7,0)*+{_2}*\frm{o}="A";
 (0,0)*+{_1}*\frm{o}="B";
 \ar @{-} "A";"B" <0pt>
\endxy} \Ea,
\Ea
\Eeq
with the representation (\ref{3: rho from RGra to differential CycW})
one recovers a wonderful result \cite{Ch} that, given any dg vector space $(W,\p)$ with the scalar product $\Theta$ as above, the associated complex $Cyc(W)$ of cyclic words in $W$ is a dg involutive Lie bialgebra with, e.g., the Lie bracket given explicitly by
$$
[(w_1\ot...\ot w_n)^{\Z^*_n}, (v_1\ot ...\ot v_m)^{\Z^*_n}]:=\hspace{100mm}
$$
$$
\hspace{20mm} \sum_{i=1}^n\sum_{j=1}^m
 \pm
\Theta(w_i,v_j) (w_1\ot ...\ot  w_{i-1}\ot v_{j+1}\ot ... \ot v_m\ot v_1\ot ... \ot v_{j-1}\ot w_{i+1}\ot\ldots\ot w_n)^{\Z^*_{n+m-2}}
$$
 %This fact plays a very important role in string topology \cite{CS1,CS2}.
 The quotient vector space,
 $$
 \overline{Cyc}(W):= Cyc(W)/\K\simeq  \bigoplus_{k\geq 1} \left(\ot^k W\right)_{\Z^*_k},
 $$
is also a $\LoBd$-algebra. On the other hand, the subspace
$$
Cyc^{\geq 3}(W):=\bigoplus_{k\geq 3} \left(\ot^k W\right)_{\Z^*_k},
$$
is a dg Lie subalgebra of $Cyc(W)$.

\subsection{Chain gravity properad and cyclic $\cA ss_\infty$ algebras}\label{3: subsec on cyclic A_infty} Let $(V,\p)$ be a dg vector space\footnote{We work in this section in the category of possibly infinite-dimensional vector spaces $V$ which are direct limits, $\displaystyle V=\lim_{\lon} V_p$, of finite-dimensional ones. Their dual vector spaces are defined as projective limits, $\displaystyle V^*=\lim_{\longleftarrow} V^*_p$.} equipped with a degree $-d$ non-degenerate scalar product
$$
\Ba{rccc}
\langle\ ,\ \rangle: & V\odot V & \lon & \K[-d]\\
                     &v_1\odot v_2     & \lon & \langle v_1 ,v_2 \rangle
\Ea
$$
which is compatible with the differential $\p$.
Hence there is an isomorphism of complexes
$$
V\lon V^*[-d]
$$
so that the vector space $V^*[-1]\simeq V[d-1]$ comes equipped with a a degree  $d-2$ pairing

$$
\Ba{rccc}
\Theta: & V^*[-1]\ot V^*[-1] & \lon & \K[d-2]]=\K[1-(3-d)]] \\
  &  \fs^{d-1}v_1\ot \fs^{d-1}v_2 & \lon & (-1)^{(d-1)|v_1|}\fs^{2d-2}\langle v_1 ,v_2 \rangle
\Ea
$$
such that
\Beqrn
\Theta(\fs^{d-1}v_2\ot \fs^{d-1}v_1) & = & (-1)^{(d-1)|v_2|}\fs^{2d-2}\langle v_2 ,v_1 \rangle  \\
  & = & (-1)^{(d-1)|v_2| +|v_1||v_2|}\fs^{2d-2}\langle v_1 ,v_2 \rangle  \\
  & = & (-1)^{(d-1)(|v_1|+|v_2|) +|v_1||v_2|}\fs^{2d-2}\Theta(\fs^{d-1}v_1, \fs^{d-1}v_2) \\
&=& (-1)^{3-d+|\fs^{d-1}v_1||\fs^{d-1}v_1|}\Theta(\fs^{d-1}v_1, \fs^{d-1}v_2).
\Eeqrn
Here we used the fact that %$\langle v_1 ,v_2 \rangle\neq 0$ only if $|v_1|+|v_2|=d$ so that
$$
(-1)^{3-d+|\fs^{d-1}v_1||\fs^{d-1}v_1|}=(-1)^{3-d+(d-1+|v_1|)(d-1+|v_1|)}=(-1)^{(d-1)(|v_1|+|v_2|) +|v_1||v_2|}
$$
Therefore, Proposition {\ref{3: Prop on repr of RGra_2}} says that the properad $\RGra_{3-d}$ acts canonically on the graded vector space $\overline{Cyc}(V^*[-1])$; in particular this space is a
$\LoB_{3-d}$-algebra and its subspace  ${Cyc}^{\geq 3}(V^*[-1])\subset \overline{Cyc}(V^*[-1])$ a dg $\Lie_{3-d}$-subalgebra with Lie bracket $[\ ,\ ]$ of degree $d-2$.

\subsubsection{\bf Lemma}\label{3: Fact on cyclic A_infty and MC elements} {\em There is a one-to-one correspondence between cyclic $\cA ss_\infty$ structures in $(V, \langle\ ,\ \rangle)$ and Maurer-Cartan elements of the dg Lie algebra   $({Cyc}^{\geq 3}(V^*[-1]), [\ ,\ ],\p)$.}

\begin{proof} This is a well-known and easy observation; hence we show only a sketch.
 A degree $3-d$ element $\ga\in {Cyc}^{\geq 3}(V^*[-1])$ gives naturally rise to a degree 1 element $\ga'$ in $\Hom(\ot^{\geq 2} V[1], V[1])$. Then the Maurer-Cartan equation
$$
\p \ga + \frac{1}{2}[\ga,\ga]=0
$$
says that sum
$$
\p + \ga'\in \Hom\left(\ot^{\geq 1} V[1], V[1]\right)
$$
is a codifferential of the tensor coalgebra $\ot^{\geq 1}(V[1])$, i.e.\ an $\cA ss_\infty$ structure in $V$. If $\ga\in (\ot ^3 V^*[-1])_{\Z_3}$, then it makes $V$ into dg associative algebra.
\end{proof}

Given a cyclic $\cA ss_\infty$ structure in a dg vector space $V$, the associated Maurer-Cartan element $\ga\in {Cyc}^{\geq 3}(V^*[-1])$  makes the vector space of ``cyclic words" $\overline{Cyc}(V^*[-1])$ into a complex with the differential
$$
\p_\ga:=\p + [\ga,\ ]
$$
It is called the cyclic Hochschild complex of the cyclic $\cA ss_\infty$-algebra $V$. The general property of T.\ Willwacher's twisting endofunctor $\tw$ (see \S 2.1) implies the following result.

\subsubsection{\bf Proposition}\label{3: Prop on action of ChGrav on Cyc} {\em The chain gravity properad $\cC h\Grav_{3-d}$ admits a canonical representation,
$$
\rho_\ga: \left(\cC h\Grav_{3-d}, \delta\right) \lon \left(\cE nd_{\overline{Cyc}(V^*[-1])}, \p_\ga\right)
$$
in the cyclic Hochschild complex of any cyclic $\cA ss_\infty$ structure in a dg vector space $V$ equipped with a non-degenerate scalar product $\odot^2 V \rar \K[-d]$.}

\sip

We are interested in this paper in applications of this general statement in string topology
of connected and simply connected closed manifolds $M$ via associated Poincar\'e models.

%%%%%%%%%%%%%%%%%%%%%%%%%%%%%%%%%%%%%%%%%%%%%%%%%%%%%%%%%
\subsection{An application to Poincar\'e duality algebras}
Let $(A,\p)$ be a finite-dimensional non-negatively graded complex equipped with a degree zero
multiplication
$$
\Ba{ccc}
 A\odot A & \lon & A\\
  a\ot b  & \lon & a\cdot b
\Ea
$$
making $A$ into a dg commutative associative algebra. It is called a {\em dg Poincar\'e duality algebra  of degree $d$}\, if it comes  equipped with a degree $-d$
 {\em orientation map}\,
 $$
 \fo: A\rar \K[-d]
 $$
 such that the induced graded symmetric scalar product,
 \Beq\label{6: scalar product in A}
\Ba{rccc}
\langle\ ,\ \rangle: & A\odot A & \lon & \K[-d]\\
                     &a\odot b     & \lon & \fo(a\cdot b)
\Ea
\Eeq
is non-degenerate and, moreover, $\fo(\p a)=0$ for any $a\in A$. We also assume that $A$ is connected and augmented,
$$
A=\K\oplus \bar{A}
$$
 with $\bar{A}$ being a positively graded vector space. It is shown in \cite{LS} that any connected and simply connected closed $d$-dimensional manifold $M$ admits a Poincar\'e duality model $A$ in degree $d$ which is quasi-isomorphic to the de Rham algebra $\Omega_M^\bu$ of $M$ as a homotopy commutative associative algebra.

\sip

The scalar product in $A$ induces an isomorphism
$$
\iota: A\lon A^*[-d]:=\Hom(A,\K)[-d]
$$
such that for any $a,b\in A$ a one has $\imath(a)(b)=\fo(a\cdot b)$. This isomorphism combined with the dualization of the multiplication map induces in turn a degree $d$ diagonal on $A$,
%$$
%\Delta: A \lon  A^*[-n] \lon  (A^* \o A^*)[-n]\lon (A[n]\ot A[n])[-n]=A\odot A [n]
%$$
$$
\Ba{rccc}
\Delta: & A & \lon & A\odot A [d]\\
    & a & \lon & \Delta(a)=: \sum a' \ot a''
\Ea
$$
which satisfies
$$
\Delta (a\cdot b)= \sum (a\cdot b') \ot b''=\sum (-1)^{|a||b'|} b'\ot (a\cdot b'') =\sum (-1)^{|a''||b|} (a'\cdot b) \ot a''= \sum  a'\ot (a''\cdot b)
$$
for any $a,b\in A$, and hence makes $A$ into a {\em Frobenius algebra}.
% In particular, we have the equalities (presented up to the obvious Koszul sign as in the formulae just %above)
%\Beq\label{6: o(abc)}
%\sum \fo(a\cdot b\cdot c') c''= \pm \fo(a\cdot b\cdot c'') c'= \pm \fo(a'\cdot b\cdot c) a''=  \pm %\fo(a''\cdot b\cdot c) a'=\pm \fo(a\cdot b'\cdot c) b''=\pm \fo(a\cdot b''\cdot c) b'.
%\Eeq
%for any $a,b,c\in A$, which we shall use later.

\sip

 As we discussed in \S {\ref{3: subsec on cyclic A_infty}} above, the graded vector space  $\overline{Cyc}({A}^*[-1])$  carries canonically a representation of the properad $\RGra_{3-d}$; in particular, it is a dg $\Lie_{3-d}$-algebra which, in accordance with Lemma {\ref{3: Fact on cyclic A_infty and MC elements}}, admits a Maurer-Cartan element
  $\ga\in (\ot ^3 {A}^*[-1])_{\Z_3}$ encoding the graded commutative associative algebra structure in ${A}$. It can be given explicitly as follows  \cite{CFL,NW},
 \Beq\label{3: gamma explicit}
  \ga=\sum_{\al,\be\ga=1}^{\dim {A}} \fo(e_\al\cdot e_\be\cdot e_\ga) s^{-1}e^\al\ot s^{-1}e^\be \ot s^{-1}e^\ga
  \Eeq
where $\{e_\al\}_{\al\in [\dim A]}$  is an arbitrary basis of ${A}$ and $e^{\al}$ the associated dual basis of $A^*$, If we choose  $e_1=1\in \K$ (resp.,  $e^1=: 1^*\in A^*$) to be the basis vector of a basis of $A$ respecting the direct sum decomposition $\K\oplus \bar{A}$  of $A$ and set $\om^*:=\fs^{d}\iota(1)\in \bar{A}^*$, then
%If $\bar{e}_\al$ is a basis of $\bar{A}$
%with $\om^*$ denoting the generator of $(A^*)^{-d}$, then
\Beqrn
\ga &=&%(s^{-1}1^*\ot s^{-1}1^* \ot s^{-1} \om^*)^{\Z^*_3} +
\sum_{\al,\be=1}^{\dim A} \fo\left({e}_\al\cdot {e}_\be\right)
(\fs^{-1}1^* \ot \fs^{-1}{e}^\al \ot \fs^{-1}{e}^\be)^{\Z^*_3}
  + \sum_{\al,\be,\ga=2}^{\dim {A}} \fo({e}_\al\cdot {e}_\be\cdot {e}_\ga) s^{-1}{e}^\al\ot s^{-1}{e}^\be \ot s^{-1}{e}^\ga
\Eeqrn
% Note that in the Lie bracket {ga,ga} only 1^* and w^* can be paired so it is enough to check
% {ga,ga}=0 only on (1\ot 1\ot w) summands --- there two pairing of w with 1 which come with opposite %signs. Next d\ga=0 because d in A is compatible with both multiplication and scalar product.
% Hence d\ga + 1/2{ga,ga}=0 indeed.

\sip

By Proposition {\ref{3: Prop on action of ChGrav on Cyc}}, {\it the chain gravity properad $\cC h\Grav_{3-d}$
acts canonically on the cyclic Hochschild complex $(\overline{Cyc}(A^*[1]), \p_\ga)$ of any Poincar\'e duality algebra $A$ in degree $d$.}

\sip

It is well-known (and easy to check) that the subspace $\overline{Cyc}(\bar{A}^*[1])\subset \overline{Cyc}(A^*[1])$ is respected by the differential $\p_\ga$ giving us the so called
{\em reduced}\, cyclic Hochschild complex of $A$. Moreover, as all black vertices of ribbon graphs from
$\cC h\Grav_{3-d}$ are at least trivalent, the operations $\rho(\Ga)$ for any $\Ga\in \cC h\Grav_{3-d}$ can not create cyclic words with letter $1^*$ when applied to elements of $\overline{Cyc}(\bar{A}^*[1])$. Hence we obtain finally the following

\subsubsection{\bf Theorem}\label{3: Theorem on action of ChGrav on Cyc(bar(A))} {\em The chain gravity properad $\cC h\Grav_{3-d}$ acts canonically and, in general, non-trivially
 on the reduced Hochschild complex $(\overline{Cyc}(\bar{A}^*[1]), \p_\ga)$
of any Poincar\'e duality algebra $A$ in degree $d$.}

\mip

More details on the non-triviality of this action will be given in the next section where we show that this action factors through a new dg {\em string topology properad}\, which has rich cohomology and which explains in a unified way  many well-known  universal operations in string topology introduced and studied  earlier in \cite{CS1,CS2,CEG,We,NW,Me2}.

%%%%%%%%%%%%%%%%%%%%%%%%%%%%%%%%%%%%%%%%%%%%%%%%%%%%%%%%%%%%%%%%%%%%%%%%%%%%%
\subsection{Gravity and equivariant cohomology of free loop spaces}
Let $LM$ stand for the space of free loops in a closed oriented $n$-dimensional manifold $M$, and
 $\bar{C}_\bu^{S^1}(LM)$) for the reduced equivariant chain complex of $LM$ with respect to the obvious $S^1$ action on $LM$; note that we work in the cohomological setting everywhere so that
 $\bar{C}_\bu^{S^1}(LM)$ is non-positively graded with the boundary differential  of cohomological degree $+1$.
 Let  $\bar{H}_\bu^{S^1}(LM)$  stand for its homology.

\sip

Assume  $A$ is a Poincar\'e model of $M$, then there is a morphism of cohomology groups
\Beq\label{3: H_S^1 into Cyc}
\bar{H}_\bu^{S^1}(LM) \lon H^{\bu}( \overline{Cyc}(\bar{A}^*[-1]), \p_\ga),
% \ \ \
%H^{\bu}( \overline{Cyc}(\bar{A}[1]), d_H^*) \lon
%\bar{H}^\bu_{S^1}(LM)
\Eeq
which is an isomorphisms if $M$ is a closed, connected and  simply connected closed manifold. This result was first proven in \cite{CEG}; another proof can be found in \cite{NW}. Therefore this isomorphism combined with Theorem {\ref{3: Theorem on action of ChGrav on Cyc(bar(A))}} implies immediately the following

\subsubsection{\bf Corollary} {\em The gravity properad $\Grav_{3-d}$ admits a representation on the
 reduced equivariant equivariant homology $\bar{H}_\bu^{S^1}(LM)$ of the free loop space $LM$ of any
 connected and simply-connected closed manifold $M$.}

\sip

We shall discuss other implications of Theorem {\ref{3: Theorem on action of ChGrav on Cyc(bar(A))}}  as well as non-trivially  of the action of  $\Grav_{3-d}$ on $\bar{H}_\bu^{S^1}(LM)$ for generic manifolds in the next section.

\bip

%%%%%%%%%%%%%%%%%%%%%%%%%%%%%%%%%%%%%%%%%%%%%%%%%%%%%%%%%%%%%%%%%%%%%%%%%%%%
%%%%%%%%%%%%%%%%%%%%%%%%%%%%%%%%%%%%%%%%%%%%%%%%%%%%%%%%%%%%%%%%%%%%%%%%%%%%

{\Large
\section{\bf String topology properad}
}

\mip

\subsection{On the joint kernel of the representation maps} Given a Poincar\'e duality algebra $A$ in degree $d$, by Theorem {\ref{3: Theorem on action of ChGrav on Cyc(bar(A))}} there is a canonical morphism of dg properads,
$$
\rho_A: \left(\cC h\Grav_{3-d},\delta\right) \lon \left(\cE nd_{ \overline{Cyc}(\bar{A}^*[-1])},\p_\ga\right).
$$
Consider the differential closure $\langle I,\delta I\rangle$ of the ideal $I$ in the non-differential properad $(\cC h\Grav_{3-d},0)$ which is generated by ribbon graphs $\Ga$ such that
\Bi
\item[(i)] every black vertex of $\Ga$ (if any) is at least four valent, or
\item[(ii)] $\Ga$ contains a boundary made of black vertices only.
\Ei
It is worth noting that the subspace of $(\cC h\Grav_{3-d},0)$ generated by graphs with property (ii) only  do {\em not}\, form an ideal as their properadic compositions with other graphs can create graphs $\Ga'$ with no boundaries made of black vertices solely; however such graphs $\Ga'$ would contain at least one black vertex of valency $\geq 4$ so that the subspace $I\subset (\cC h\Grav_{3-d},0)$ generated by graphs of {\em both}\, types (i) and (ii) is a properadic ideal indeed.

\sip

It is clear that
$$
\langle I,\delta I\rangle \in \Ker \rho_A
$$
for any Poincar\'e duality algebra $A$ because
\Bi
 \item[$(i)$] graphs $\Ga$ from (i) above do not contribute into $\rho_A$ as the Maurer-Cartan element $\ga\in \overline{Cyc}({A}^*[-1])$ (which decorates black vertices upon the representation)
is a linear combination of cyclic words with precisely three letters from ${A}^*[-1]$ (see (\ref{3: gamma explicit})),
\item[$(ii)$] graphs $\Ga$ from (ii) above do not contribute since we work with the quotient spaces, $\overline{Cyc}(V)=Cyc(V)/\K$, which do not contain the ``empty" cyclic word, i.e.\ the one with no letters from a vector space $V$.
\Ei

These observations together with Theorem {\ref{3: Theorem on action of ChGrav on Cyc(bar(A))}} and  isomorphism (\ref{3: H_S^1 into Cyc}) imply the following

\subsubsection{\bf Definition-Theorem}\label{4: def-ptop on ST} {\em The quotient dg properad
$$
\ST_{d}:= \cC h\Grav_{d}/ \langle I,\delta I\rangle, \ \ \ d\in \Z.
$$
is called the {\em chain string topology}\, properad. The properad $\ST_{3-d}$ admits a canonical representation on the reduced cyclic Hochschild complex $(\overline{Cyc}(\bar{A}^*[-1]),\p_\ga)$ of a Poincar\'e duality algebra $A$ in degree $d$. In particular, its cohomology properad $H^\bu(\ST_{3-d})$ (called the {\em string topology} properad) acts on the
 reduced equivariant equivariant homology $\bar{H}_\bu^{S^1}(LM)$ of the free loop space $LM$ of any
 connected and simply-connected closed manifold $M$.}

\mip

The induced differential $\delta$ in $\ST_d$ acts only on white vertices splitting them as follows
$$
\Ba{c}\resizebox{5mm}{!}{  \xy
 (0,0)*+{_i}*\frm{o}="B";
\endxy} \Ea
\lon
\Ba{c}\resizebox{5mm}{!}{  \xy
 (0,6)*{\bu}="A";
 (0,0)*+{_i}*\frm{o}="B";
 \ar @{-} "A";"B" <0pt>
\endxy} \Ea,
$$
and redistributing edges among the new vertices in such a way that the newly created black vertex becomes strictly trivalent. In particular, the differential in $\ST_d$ acts trivially on univalent white vertices. Let us check next the non-triviality of the cohomology properad $H^\bu(\ST_d)$.

\subsection{$H^\bu(\ST_d)$ and involutive Lie bialgebras} The r.h.s.\ in the formula (\ref{2: delta on ribbon CoLie corolla}) for
$
\delta \xy
 (0.5,1)*{^{{^{\bar{1}}}}},
(0.5,5)*{^{{^{\bar{2}}}}},
 (0,-2)*+{_{_1}}*\frm{o}="A";
% (0,0)*+{_1}*\frm{o}="B";
%(0,-2)*{\circ}="A";
%(0,-2)*{\circ}="B";
"A"; "A" **\crv{(7,7) & (-7,7)};
\endxy
$ in $\cC h\Grav_d$ belongs to the ideal $\langle I,\delta I\rangle$ and hence vanishes in $\ST_d$. Therefore there is a morphism of dg properads (cf.\ (\ref{3: morphism from LoB to RGra}))
$$
\Ba{rccc}
i: & (\LoB_{d},0) & \lon & (\ST_{d},\delta)\vspace{2mm}\\
&
\Ba{c}\begin{xy}
 <0mm,-0.55mm>*{};<0mm,-2.5mm>*{}**@{-},
 <0.5mm,0.5mm>*{};<2.2mm,2.5mm>*{}**@{-},
 <-0.48mm,0.48mm>*{};<-2.2mm,2.5mm>*{}**@{-},
 <0mm,0mm>*{\bu};<0mm,0mm>*{}**@{},
 <0mm,-0.55mm>*{};<0mm,-3.8mm>*{_{_1}}**@{},
 <0.5mm,0.5mm>*{};<2.7mm,3.1mm>*{^{^{\bar{2}}}}**@{},
 <-0.48mm,0.48mm>*{};<-2.7mm,3.1mm>*{^{^{\bar{1}}}}**@{},
 \end{xy}\Ea
 &\lon &
 \Ba{c}\resizebox{7mm}{!}{ \xy
 (0.5,1)*{^{{^{\bar{1}}}}},
(0.5,5)*{^{{^{\bar{2}}}}},
 (0,-2)*+{_{_1}}*\frm{o}="A";
"A"; "A" **\crv{(7,7) & (-7,7)};
\endxy}\Ea
\vspace{2mm}\\
&
\Ba{c}\begin{xy}
 <0mm,0.66mm>*{};<0mm,3mm>*{}**@{-},
 <0.39mm,-0.39mm>*{};<2.2mm,-2.2mm>*{}**@{-},
 <-0.35mm,-0.35mm>*{};<-2.2mm,-2.2mm>*{}**@{-},
 <0mm,0mm>*{\bu};<0mm,4.1mm>*{^{^{\bar{1}}}}**@{},
   <0.39mm,-0.39mm>*{};<2.9mm,-4mm>*{^{_2}}**@{},
   <-0.35mm,-0.35mm>*{};<-2.8mm,-4mm>*{^{_1}}**@{},
\end{xy}\Ea
 & \lon &
 \Ba{c}\resizebox{10mm}{!}{  \xy
 (3.5,4)*{^{\bar{1}}};
 (7,0)*+{_2}*\frm{o}="A";
 (0,0)*+{_1}*\frm{o}="B";
 \ar @{-} "A";"B" <0pt>
\endxy} \Ea,
\Ea
$$
implying that the reduced cyclic Hochschild complex (and hence its cohomology) is always an involutive Lie bialgebra. In particular, this morphism  and Theorem {\ref{4: def-ptop on ST}} imply a purely algebraic construction of the involutive Lie bialgebra $\LoB_{3-d}$-structure   on $\bar{H}_\bu^{S^1}(LM)$; this construction was first found by  X.\ Chen, F.\ Eshmatov, and W.\ L.\ Gan in \cite{CEG}. It was proven by F.\ Naef and T.\ Willwacher in \cite{NW} that this algebraic construction of $\LoB_{3-d}$ structure on $\bar{H}_\bu^{S^1}(LM)$ agrees precisely with the original geometric construction of M.\ Chas and D.\ Sullivan in \cite{CS1,CS2}.

\subsection{An injection of the gravity operad into $H^\bu(\ST_d)$}
 The chain gravity properad
$\cC h\Grav_d$ contains a sub-operad $\tw\cR \cT ree_d$ spanned by genus zero ribbon graphs with precisely one boundary, i.e.\ by ribbon trees. It is proven in \cite{Wa} that its cohomology operad
$$
H^\bu(\tw\cR\cT ree_d)=grav_v
$$
can be identified with the E.\ Getzler's gravity operad (\ref{1: grav_d operad}). Moreover, every cohomology class in $H^\bu(\cR\cT ree_d)$ admits a representation in terms of a linear combination of ribbon trees with all white vertices univalent and black vertices trivalent as, for example,
the following ones
$$
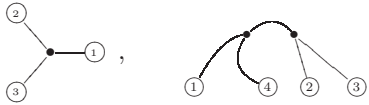
\Ba{c}\resizebox{14mm}{!}{
\mbox{$\xy
 (0,0)*{\bu}="C";
  (7.9,0)*{_{_1}}*+\frm{o}="1";
(-6,7)*{_{_2}}*+\frm{o}="2";
(-6,-7)*{_{_3}}*+\frm{o}="3";
 \ar @{-} "1";"C" <0pt>
 \ar @{-} "2";"C" <0pt>
  \ar @{-} "3";"C" <0pt>
\endxy$}}
\Ea
, \ \ \ \
\Ba{c}\resizebox{25mm}{!}{ \xy
(-1.5,5)*{}="1",
(1.5,5)*{}="2",
(9,5)*{}="3",
 (0,0)*{\bu}="A";
  (9,0)*{\bu}="O";
 (-10,-10)*+{_1}*\frm{o}="B";
 (4,-10)*+{_4}*\frm{o}="E";
  (12,-10)*+{_2}*\frm{o}="C";
   (21,-10)*+{_3}*\frm{o}="D";
 "A"; "B" **\crv{(-5,-0)}; % (5,5) = tangent point
 %"A"; "B" **\crv{(-5,-6)};
  %"A"; "D" **\crv{(5,-0.5)};
   %"A"; "B" **\crv{(5,-1)};
  "A"; "E" **\crv{(-5,-7)};
   "A"; "O" **\crv{(5,5)};
\ar @{-} "O";"C" <0pt>
\ar @{-} "O";"D" <0pt>
 \endxy}
 \Ea
$$
Therefore we can make the following conclusion.

\subsubsection{\bf Lemma} {\em The gravity operad $grav_d$ injects into the string topology properad $H^\bu(\ST_d)$}.

\mip

This Lemma together with Theorem {\ref{4: def-ptop on ST}} imply an action of the gravity properad $grav_{3-d}$ on the
 reduced equivariant equivariant homology $\bar{H}_\bu^{S^1}(LM)$ of the free loop space $LM$ of any
 connected and simply-connected closed manifold $M$. This is a purely algebraic counterpart of the geometric construction by C.\ Westerland in \cite{We}. Conjecturally, both constructions describe identical string topology operations.

\subsection{$H^\bu(\ST_d)$ and homotopy involutive Lie bialgebras} Consider the following four graphs in $\ST_{d}$, $\forall d\in \Z$,
\Beq\label{4: four trinary cycles in ST}
\xy
(0,4)*{\bu}="1";
(-6,-4)*{\ }*+\frm{o}="2";
%(-5,-4)*{\circ}="2";
(0,-4)*{\ }*+\frm{o}="3";
%(0,-4)*{\circ}="3";
(6,-4)*{\ }*+\frm{o}="4";
%(5,-4)*{\circ}="4";
\ar @{-} "1";"2" <0pt>
\ar @{-} "1";"3" <0pt>
\ar @{-} "1";"4" <0pt>
\endxy\in \ST_{d}(1,3),  \
\ \ \ \xy
(0,5)*{\bu}="1";
(0,-4)*{\ }*+\frm{o}="3";
%(0,-4)*{\circ}="3";
"1";"3" **\crv{(4,0) & (4,1)};
"1";"3" **\crv{(-4,0) & (-4,-1)};
\ar @{-} "1";"3" <0pt>
\endxy\in \ST_{d}(3,1),  \ \
\xy
(0,5)*{\bu}="1";
(-5,-4)*{\ }*+\frm{o}="3";
%(-5,-4)*{\circ}="3";
(5,-4)*{\ }*+\frm{o}="4";
%(5,-4)*{\circ}="4";
"1";"3" **\crv{(-3,5) & (5,4)};
"1";"3" **\crv{(-5,2) & (-5,2)};
\ar @{-} "1";"4" <0pt>
\endxy\in \ST_{d}(2,2),
\ \ \ \xy
(0,5)*{\bu}="1";
(0,-4)*{\ }*+\frm{o}="3";
%(0,-4)*{\circ}="3";
"1";"3" **\crv{(-5,2) & (5,2)};
"1";"3" **\crv{(5,2) & (-5,2)};
"1";"3" **\crv{(-7,7) & (-7,-7)};
%\ar @{-} "1";"3" <0pt>
\endxy\in \ST_{d}(1,1),
%
%
%
%\ \ \ \xy
%(0,5)*{\ast}="1";
%(0,-5)*{\circ}="3";
%"1";"3" **\crv{(-5,2) & (5,2)};
%"1";"3" **\crv{(5,2) & (-5,2)};
%"1";"3" **\crv{(-7,0) & (-7,-1)};
%\ar @{-} "1";"3" <0pt>
%\endxy
%\ \ \ \xy
%(0,5)*{\ast}="1";
%(0,-5)*{\circ}="3";
%
%"1";"3" **\crv{(-3,5) & (5,4)};
%"1";"3" **\crv{(5,2) & (-5,2)};
%\ar @{-} "1";"3" <0pt>
%\endxy
\Eeq
whose white vertices and boundaries are symmetrized (resp., skewsymmetrized) for $d$ odd (resp., $d$ even) so that we can omit their numerical labels in pictures.
% asterisk has degree 1+D, half-edges have degree -D. Set 1+D=d, so the "black" vertex has degree d while "half" edges have degree
%1-d as expected in ST_d.
These elements are obviously cycles in $\ST_{d}$ which --- due to the (skew)symmetrizations --- can not be coboundaries. Hence they give us four non-trivial cohomology classes in $H^\bu(\ST_{d})$
which, for $d\leq 1$,  give us in turn four universal string topology operations on $\bar{H}_\bu^{S^1}(LM)$ for any
closed, connected and simply-connected $(3-d)$-dimensional manifold $M$. These four operations have been discovered in \cite{Me2} using a completely different technical gadget called the properad of ribbon {\em hyper}graphs. Moreover, the first half of Proposition 6.3.1 from \cite{Me2} together with with the above observation imply the following result.

\subsubsection{\bf Proposition} {\em There is a morphism of dg properads
$$
j: (\HoLoB_{d-1}, \delta) \lon (\ST_{d},\delta)
$$
which vanishes on all generators (\ref{2: generators of HoLoB}) of $\HoLoB_{d-1}$ except the following quartette,
$$
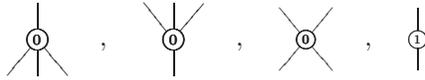
\Ba{c}\resizebox{8mm}{!}{\xy
(-5,-6)*{};
(0,0)*+{_0}*\cir{}
**\dir{-};
(0,-6)*{};
(0,0)*+{_0}*\cir{}
**\dir{-};
(5,-6)*{};
(0,0)*+{_0}*\cir{}
**\dir{-};
(0,6)*{};
(0,0)*+{_0}*\cir{}
**\dir{-};
%
%(0,8)*{_1};
%(-5,-8)*{_1};
%(0,-8)*{_2};
%(5,-8)*{_3};
\endxy}\Ea\ \ , \ \
\Ba{c}\resizebox{8mm}{!}{\xy
(-5,6)*{};
(0,0)*+{_0}*\cir{}
**\dir{-};
(0,6)*{};
(0,0)*+{_0}*\cir{}
**\dir{-};
(5,6)*{};
(0,0)*+{_0}*\cir{}
**\dir{-};
(0,-6)*{};
(0,0)*+{_0}*\cir{}
**\dir{-};
%
%(0,-8)*{_1};
%(-5,8)*{_1};
%(0,8)*{_2};
%(5,8)*{_3};
\endxy}\Ea\ \ , \ \
\Ba{c}\resizebox{7mm}{!}{\xy
(-5,-6)*{};
(0,0)*+{_0}*\cir{}
**\dir{-};
(5,6)*{};
(0,0)*+{_0}*\cir{}
**\dir{-};
(5,-6)*{};
(0,0)*+{_0}*\cir{}
**\dir{-};
(-5,6)*{};
(0,0)*+{_0}*\cir{}
**\dir{-};
%
%(-5,8)*{_1};
%(5,8)*{_2};
%(-5,-8)*{_1};
%(5,-8)*{_2};
\endxy}\Ea\ \ , \ \
\Ba{c}\resizebox{2.5mm}{!}{\xy
(0,-6)*{};
(0,0)*+{_1}*\cir{}
**\dir{-};
(0,6)*{};
(0,0)*+{_1}*\cir{}
**\dir{-};
%
%(0,8)*{_1};
%(0,-8)*{_1};
\endxy}\Ea
$$
which is sent, respectively, by $j$ into the list (\ref{4: four trinary cycles in ST}) of cycles in $\ST_{d}$.}
% Generators of HoLoBd have degree 1-d(m+n+2a-2).
% hence these four geberatots of HoLoB_{d-1} have degree 1-(d-1)(4-2)=3-2d,
% which agrees with degrees d+3(1-d)=3-2d of the graphs with one black vertex and three edges!
% Everything is fine.

\mip

This Proposition together with Theorem {\ref{4: def-ptop on ST}} gives us a new surprising proof of the second part of the Proposition 6.3.1  in \cite{Me2} asserting a non-trivial role of $\HoLoB_{2-d}$ in string topology of generic $d$-dimensional manifolds.

\bip

There are many more non-trivial cohomology classes in $H^\bu(\ST_d)$ (and hence universal operations in string topology), e.g., the following ones
$$
%\xy
%(0,5)*{\bu}="1";
%(-5,-4)*{\ }*+\frm{o}="3";
%(10,-4)*{\ }*+\frm{o}="5";
%(7,3)*{\bu}="1'";
%"1";"3" **\crv{(-3,5) & (5,4)};
%"1";"3" **\crv{(-5,2) & (-5,2)};
%\ar @{-} "1";"1'" <0pt>
%\ar @{-} "3";"1'" <0pt>
%\ar @{-} "5";"1'" <0pt>
%\endxy
%\ \ ,
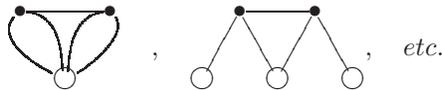
\xy
(-6,5)*{\bu}="1";
(0,-4)*{\ }*+\frm{o}="3";
(6,5)*{\bu}="2";
%(15,-4)*{\ }*+\frm{o}="4";
%(5,-4)*{\circ}="4";
"1";"3" **\crv{(1,5) & (-1,-4)};
"1";"3" **\crv{(-10,2) & (-5,-2)};
"2";"3" **\crv{(10,2) & (1,-4)};
"2";"3" **\crv{(-1,5) & (1,-5)};
\ar @{-} "1";"2" <0pt>
\endxy
\ \ \ , \ \ \
\xy
(0,5)*{\bu}="1";
(-5,-4)*{\ }*+\frm{o}="3";
%(-5,-4)*{\circ}="3";
(5,-4)*{\ }*+\frm{o}="4";
(15,-4)*{\ }*+\frm{o}="5";
(10,5)*{\bu}="1'";
%"1";"3" **\crv{(-3,5) & (5,4)};
%"1";"3" **\crv{(-5,2) & (-5,2)};
\ar @{-} "1";"4" <0pt>
\ar @{-} "1";"1'" <0pt>
\ar @{-} "4";"1'" <0pt>
\ar @{-} "5";"1'" <0pt>
\ar @{-} "1";"3" <0pt>
\endxy
,\ \ \ etc.
$$
where the first (resp. second) ribbon graph has all boundaries (resp.\ vertices) skew-symmetrized for $d$ odd/symmetrized for $d$ even.
We conclude that the string topology properad $H^\bu(\ST_d)$ is highly non-trivial and is worth of further study. It is also an interesting open problem to study the cohomology group of the deformation complex of the morphism $i$  discussed above,
$$
\Def(\LoB_{d} \stackrel{i}{\rar} \ST_d), %\ \ \ \ \Def(\HoLoB_{d-1}\lon \ST_d)
$$
as it comes equipped with a  morphism from the cohomology group (cf.\ \cite{MW2})
$$
H^\bu(\GC_{2d}) \lon H^\bu\left(\Def(\LoB_{d} \stackrel{i}{\rar} \ST_d)\right),
 %\ \ \ H^\bu(\GC_{2(d-1)}) \lon\Def(\HoLoB_{d-1}\stackrel{j}{\rar} \ST_d)
$$
of the mysterious even Kontsevich graph complex $\GC_{2d}$ which have been introduced in \cite{Ko} and studied in \cite{W}.

\bip

\bip

%%%%%%%%%%%%%%%%%%%%%%%%%%%%%%%%%%%%%%%%%%%%%%%%%%%%%%%%%%%%%%%%%%%%%%%%

\def\cprime{$'$}

\end{document}